\documentclass[12pt]{article}
\usepackage{latexsym}
\usepackage{graphicx}
\usepackage{cite}
\usepackage{enumitem}
\usepackage{amssymb,mathrsfs,amsmath}
\usepackage{booktabs}
\usepackage[justification=centering]{caption}
\usepackage{subfigure}
\usepackage{float}
\usepackage{color}
\usepackage{arydshln}

\usepackage{mathtools,amsthm}
\usepackage[normalem]{ulem}
\usepackage{mathabx}
\usepackage{acronym}

\usepackage[margin=1.5in]{geometry}%
\newtheorem{theorem}{Theorem}[part]
\newtheorem{definition}{Definition}[part]
\newtheorem{proposition}{Proposition}[part]

\newtheorem{lemma}{Lemma}[part]
\newtheorem{corollary}{Corollary}[part]
\newtheorem{remark}{Remark}[part]

\def\top{\mathsf{T}}

\begin{document}

\title{T-positive semidefiniteness of third-order symmetric tensors and T-semidefinite programming}

\author{Meng-Meng Zheng\thanks{School of Mathematics, Tianjin University,
Tianjin 300072, P.R. China. Email: zmm941112@tju.edu.cn.} \and
Zheng-Hai Huang\thanks{Corresponding Author. School of Mathematics,
Tianjin University, Tianjin 300072, P.R. China. Email:
huangzhenghai@tju.edu.cn. Tel:+86-22-27403615 Fax:+86-22-27403615}
\and Yong Wang\thanks{School of Mathematics, Tianjin University, Tianjin 300072, P.R. China.
Email: wang\_yong@tju.edu.cn.}}
\date{}

\maketitle

\begin{abstract}
 The T-product for third-order tensors has been used extensively in the literature. In this paper, we first introduce the first-order and second-order T-derivatives for the multi-vector real-valued function with the tensor T-product; and inspired by an equivalent characterization of a twice continuously T-differentiable multi-vector real-valued function being convex, we present a definition of the T-positive semidefiniteness of third-order symmetric tensors. After that, we extend many properties of positive semidefinite matrices to the case of third-order symmetric tensors. In particular, analogue to the widely used semidefinite programming (SDP for short), we introduce the semidefinite programming over the third-order symmetric tensor space (T-semidefinite programming or TSDP for short), and provide a way to solve the TSDP problem by converting it into an SDP problem in the complex domain. Furthermore, we give several examples which can be formulated (or relaxed) as TSDP problems, and report preliminary numerical results for two unconstrained polynomial optimization problems. Experiments show that finding the global minimums of polynomials via the TSDP relaxation outperforms the traditional SDP relaxation for the test examples.

\noindent {\bf Key words:}\hspace{2mm} T-product, T-positive semidefiniteness, T-semidefinite cone,  T-semidefinite programming, polynomial optimization.  \vspace{3mm}

\noindent {\bf Mathematics Subject
Classifications (2010):}\hspace{2mm} 15A69, 90C22 \vspace{3mm}
\end{abstract}

\section{Introduction}
\label{intro}
With the availability of inexpensive storage and advances in instrumentation, the data collected and stored now is more complex than ever before. Especially in  practical problems such as psychometrics, signal processing, computer vision, data mining, graphical analysis, neuroscience and so on, it is usually necessary to store information in a multidimensional array, and then use the multidimensional structure to compress, sort, and/or manipulate the data.
Among the many problems described by high-dimensional arrays (or tensors), third-order tensors have become increasingly prevalent in recent years with the emergence of the tensor T-product, which is a new type of multiplication between third-order tensors introduced by Kilmer, Martin, and Perrone \cite{KMP-08}. The tensor T-product has shown to be a useful tool arising in a wide variety of application areas, including, but not limited to, image processing \cite{KSB-18,KBHH-13,MSL-13,SKH-16,TM-18,ZLLZ-18}, computer vision \cite{BG-19,HKBH-13,SHKM-13,XTZL-18,YGXG-19}, signal processing, low rank tensor recovery and robust tensor PCA \cite{CYH-16,KXL-18,LCZ-18,LLC-19,LFCL-19,SHS-19}, and data completion and denoising \cite{EAH-15,HTZXY-17,HYZX-17,LAAWW-16,MG-18-2,MG-18,QJH-19,WLJ-19,YHHH-16,ZHJ-18,Z-17,Z-15,ZEAHK-14}, because the tensor T-product provides an effective approach to transform the tensor multiplication into block diagonal matrix multiplication in the discrete Fourier domain.

Since  Kilmer, Martin, and Perrone \cite{KMP-08} introduced the new type of multiplications between two third-order tensors so as to devise new types of factorizations for tensors to be easily used in applications, the exploration of the algebraic properties of T-products has been in progress. Specifically, in \cite{KM-11} some factorization strategies {were} established for third-order tensors via the tensor T-product. In \cite{B-10} and \cite{KBHH-13}, the authors provided useful frameworks to better view the action of the third-order tensors upon a set of matrices. In \cite{MQW-19}, a lot of familiar tools of linear algebra were extended to the third-order tensors, including the T-Jordan canonical form, tensor decomposition
theory, T-group inverse and T-Drazin inverse, and so on. In addition, Lund \cite{L-18} proposed the definition of tensor functions based on the T-product of third-order F-square tensors, which was found to be of great use in stable tensor neural networks for rapid deep learning \cite{NHAK-18}, and then Miao, Qi and Wei generalized the tensor T-function from F-square third order tensors to rectangular tensors in \cite{MQW-19-1}.

It is well known that the positive (semi)definite (P(S)D for short) matrix is an important class of matrices, which has a wealth of theoretical results and applications. Actually, P(S)D matrices can be used for inequality proof, eigenvalue solving, extremum solving, system stability discrimination, and so on. As a result, P(S)D matrices have been applied in various fields, such as numerical analysis, optimization theory, probability and statistics, operations research, control theory, mechanics, electricity, information science and technology, management science and engineering, and so on. More information about P(S)D matrices can refer to the monographs \cite{R-07,HJ-13}. In addition, the semidefinite programming (SDP for short) as one of the important applications of PSD matrices has received extensive attention. Especially with the appearance of some effective algorithms for SDP problems \cite{Sedumi-99,SDPNAL+-19,SDPT3-99}, SDP problems have been increasingly arisen in practical applications. There are rich and mature results for SDPs in both theory and algorithms. See \cite{H-00,T-01,VB-96} and references therein.

Motivated by that mentioned above, we extend the positive (semi)definiteness of matrices and the SDP problem to the case of third-order tensors. Our contribution is threefold:
\begin{itemize}
\item[1)] We explore the Fr\'{e}chet derivatives of the multi-vector real-valued function based on the inner product and the tensor T-product. We establish the necessary and sufficient conditions for a multi-vector real-valued function being first-order and second-order T-differentiable, respectively, and present the exact forms of the T-derivatives. In particular, we propose a second-order condition to judge the convexity of the multi-vector real-valued function under the premise that the function is twice continuously T-differentiable.
\item[2)] We give a definition of a third-order symmetric tensor being T-positive semidefinite (T-PSD for short) inspired by the second-order T-derivative, and show that the new definition is equivalent to the one given by \cite[Definition 2.7]{KBHH-13} and the one by \cite[Definition 15]{MQW-19} in real case. In particular, we show that the set of symmetric T-PSD tensors is a nonempty, closed, convex, pointed and self-dual cone, and extend many properties of PSD matrices to the case of third-order T-PSD tensors, including some results related to the T-eigenvalue decomposition, the T-roots, the T-Schur complement, and so on.
\item[3)] As an important application of the T-positive semidefiniteness of third-order symmetric tensors, we introduce the semidefinite programming over the third-order symmetric tensor space (T-semidefinite programming or TSDP for short) and show that a TSDP problem of size $m\times m\times p$ can be transformed into an SDP problem with $\frac{p+1}{2}$ or $\frac{p+2}{2}$ blocks of matrices of $m\times m$ in the complex domain. Then we present several examples which can be modeled (or relaxed) as TSDPs, such as minimizing the maximum T-eigenvalue of a third-order tensor, minimizing the spectral norm of a third-order tensor, minimizing the nuclear norm of a third-order tensor, integer quartic programming and calculating the global lower bound of a polynomial. Besides, we report preliminary numerical results for solving the unconstrained polynomial optimization problem via the TSDP relaxation, which can achieve higher accuracy and consume less time compared with the traditional SDP relaxation.
\end{itemize}

The rest of our paper is organized as follows. In Section \ref{Sect. 1}, some notation and basic results are reviewed. In Section \ref{Sect. 2}, we explore the T-derivatives for the multi-vector real-valued function and the relationship between a new type of tensors: the T-Hessian Tensor and convexity of the multi-vector real-valued function. In Section \ref{Sect. 3}, we give the definition of the symmetric T-PSD tensor, then discuss some characterizations and properties of symmetric T-PSD tensors, and investigate the set of T-PSD tensors. In Section \ref{Sect. 4}, we introduce and study the TSDP; and convert the TSDP into the corresponding SDP in the complex domain. We also present several examples for applications and report preliminary numerical results. Finally, we sum up the conclusions and do some further discussions in Section \ref{Sect. 5}.
\section{Preliminary}
\label{Sect. 1}In this section, we give some notation and basic results.

\subsection{Notation}
\label{Sect. 2.1}
Throughout this paper, we use small letters $a,b,\ldots$ for scalars, small bold letters $\mathbf{a}, \mathbf{b}, \ldots$ for vectors, capital letters $A,B,\ldots$ for sets, capital bold letters $\mathbf{A}$, $\mathbf{B}$, $\ldots$ for matrices, and calligraphic letters $\mathscr{A}, \mathscr{B}, \ldots$ for tensors. For any positive integer $n$, denote $[n] := \ \{1,2,\ldots , n\}$. Let $\mathbb{R}^n := \ \{\mathbf{x}:=(x_1,x_2,\ldots ,x_n)^\top: x_i\in \mathbb{R}\ \mbox{\rm for all}\ i\in [n] \}$ and $\mathbb{C}^n := \ \{\mathbf{x}:=(x_1,x_2,\ldots ,x_n)^\top: x_i\in \mathbb{C}\ \mbox{\rm for all}\ i\in [n] \}$ where $\mathbb{R}$ ($\mathbb{C}$) is the set of real (complex) numbers. Let $m$, $n$ and $p$ be positive integers. $\mathbb{R}^{m\times n}$ and $\mathbb{R}^{m\times n\times p}$ denote the sets consisting of all real matrices of size $m\times n$ and all real tensors of size $m\times n\times p$, respectively. Let $\mathbb{N}$ denote the set of nonnegative integers. For $\alpha\in \mathbb{N}^n$, denote $|\alpha|:=\alpha_1+\alpha_2+\cdots+\alpha_n$. For any $\mathbf{x}\in \mathbb{R}^n$ and $\alpha\in \mathbb{N}^n$, $\mathbf{x}^\alpha$ means $x_1^{\alpha_1}\cdots x_n^{\alpha_n}$, and $\mathbf{x}^\top $ represents the transpose of $\mathbf{x}$. For any $\mathscr{A}$, $\mathscr{B}\in \mathbb{R}^{m\times n\times p}$, the inner product between $\mathscr{A}$, $\mathscr{B}$ is denoted as $\mathscr{A}\bullet\mathscr{B}=\langle\mathscr{A},\mathscr{B}\rangle:=\sum_{i,j,k}a_{ijk}b_{ijk},$
and the Frobenius norm associated with the above inner product is $\|\mathscr{A}\|=\sqrt{\mathscr{A}\bullet\mathscr{A}}$. Specially, any matrix $\mathbf{A}\in\mathbb{R}^{n\times p}$ can be regarded as a tensor $\mathscr{A}\in\mathbb{R}^{n\times1\times p}$ with the $i$-th frontal slice of $\mathscr{A}$ being the $i$-th column of $\mathbf{A}$ for all $i\in [p]$.

Recall that a complex matrix $\mathbf{A}$ is said to be symmetric (Hermitian) if and only if $\mathbf{A}^\top=\mathbf{A}$ ($\mathbf{A}^H=\mathbf{A}$), where $\mathbf{A}^\top$ ($\mathbf{A}^H$) represents the transpose (conjugate transpose) of $\mathbf{A}$. We denote the set consisting of all real symmetric (complex Hermitian) matrices of size $n\times n$ as $S\mathbb{R}^{n\times n}$ ($H\mathbb{C}^{n\times n}$). For any $x\in \mathbb{C}$ and $\mathbf{X}:=(x_{ij})\in \mathbb{C}^{m\times n}$, $\overline{x}$ denotes the conjugate of $x$ and $\overline{\mathbf{X}}:=(\overline{x}_{ij})$ denotes the conjugate of the matrix $\mathbf{X}$.
Let $\mathbf{U}\succeq(\succ)0$ represent that $\mathbf{U}$ is (Hermitian) positive semidefinite (positive definite) for any $\mathbf{U}\in H\mathbb{C}^{n\times n}$, and $S\mathbb{R}^{n\times n}_{++}( S\mathbb{R}^{n\times n}_{+})$ denote the set of all real symmetric positive (semi)definite matrices of size $n\times n$, while $H\mathbb{C}^{n\times n}_{++}( H\mathbb{C}^{n\times n}_{+})$ denotes the set of all complex Hermitian positive (semi)definite matrices of size $n\times n$ . ``$\otimes$" denotes the Kronecker product between two matrices and ``$\cdot$"  means standard matrix product.

\subsection{ Tensor T-product, transpose and inverse}
\label{Sect. 2.2}
For a third-order tensor $\mathscr{A}\in \mathbb{R}^{m\times n\times p}$, a new perspective was proposed in \cite{KMP-08,KM-11} based on treating $\mathscr{A}$ as a stack of frontal slices, which were denoted as $\mathbf{A}^{(k)}\in \mathbb{R}^{m\times n}$ for all $k\in [p]$. Furthermore, several operators on $\mathscr{A}\in \mathbb{R}^{m\times n\times p}$ were introduced as follows:
\begin{eqnarray*}
bcirc(\mathscr{A}):=\left[\begin{array}{ccccc}
\mathbf{A}^{(1)} & \mathbf{A}^{(p)} & \mathbf{A}^{(p-1)} & \cdots & \mathbf{A}^{(2)} \\
\mathbf{A}^{(2)} & \mathbf{A}^{(1)} & \mathbf{A}^{(p)} & \cdots & \mathbf{A}^{(3)} \\
\vdots & \ddots & \ddots & \ddots & \vdots \\
\mathbf{A}^{(p)} & \mathbf{A}^{(p-1)} & \cdots & \mathbf{A}^{(2)} & \mathbf{A}^{(1)}
\end{array}
\right],\;\;  unfold(\mathscr{A}):=\left[\begin{array}{c}
 \mathbf{A}^{(1)}\\
 \mathbf{A}^{(2)}\\
 \vdots\\
 \mathbf{A}^{(p)}
\end{array}
\right],
\end{eqnarray*}
$fold(unfold(\mathscr{A})):=\mathscr{A}$, and $bcirc^{-1}(bcirc(\mathscr{A})):=\mathscr{A}$.

With the help of the above operators, the following definitions and properties were given in \cite{KBHH-13} (see also
\cite{HKBH-13,KMP-08,KM-11,MQW-19}).

\begin{definition}\label{T-product}\cite[Definition 2.5]{KBHH-13}
(T-product) Let $\mathscr{A}\in \mathbb{R}^{m\times n\times p}$ and $\mathscr{B}\in \mathbb{R}^{n\times s\times p}$ be two real tensors. Then the T-product $\mathscr{A}\ast\mathscr{B}$ is an $m\times s\times p$ real tensor defined by
$\mathscr{A}\ast\mathscr{B}:=fold(bcirc(\mathscr{A})unfold(\mathscr{B})).$
\end{definition}

\begin{definition}\label{transpose}\cite[Definition 2.7]{KBHH-13}
(Transpose and conjugate transpose) If $\mathscr{A}$ is a third-order tensor of size $m\times n\times p$, then the transpose $\mathscr{A}^\top$ is obtained by transposing each of the frontal slices and then reversing the order of transposed frontal slices $2$ through $p$. The conjugate transpose $\mathscr{A}^H$ is obtained by conjugate transposing each of the frontal slices then reversing the order of transposed frontal slices $2$ through $p$.
\end{definition}

For any $\mathscr{A}\in \mathbb{R}^{n\times n\times p}$, we say $\mathscr{A}$ is a symmetric tensor if and only if $\mathscr{A}^\top=\mathscr{A}$. The set consisting of all the real symmetric tensor of size $n\times n\times p$ is denoted by $S\mathbb{R}^{n\times n\times p}$.

\begin{definition}\label{inverse}\cite[Definition 2.8,2.10]{KBHH-13}
(Identity tensor and inverse) The $n\times n\times p$ identity tensor $\mathscr{I}_{nnp}$ is the tensor whose first frontal slice is the $n\times n$ identity matrix $\mathbf{I}_{n\times n}$, and whose other frontal slices are all zeroes. For a frontal square tensor $\mathscr{A}\in \mathbb{R}^{n\times n\times p}$, we say $\mathscr{A}$ is nonsingular if it has inverse tensor $\mathscr{B}(=\mathscr{A}^{-1})$, provided that $\mathscr{A}\ast\mathscr{B}=\mathscr{B}\ast\mathscr{A}=\mathscr{I}_{nnp}.$
\end{definition}

It is easy to check that $\mathscr{A}\ast\mathscr{I}_{nnp}=\mathscr{I}_{mmp}\ast\mathscr{A}=\mathscr{A}$ for any $\mathscr{A}\in \mathbb{R}^{m\times n\times p}$. In addition, it should be noticed that the invertibility of the third-order frontal square tensor $\mathscr{A}$ is equivalent to the invertibility of the matrix $bcirc(\mathscr{A})$, which can be seen from the following lemma \cite[Lemma 3]{MQW-19}.

\begin{lemma}\label{lemma-1}
Suppose that $\mathscr{A}\in \mathbb{R}^{n\times n\times p}$ and $\mathscr{B}\in \mathbb{R}^{n\times s\times p}$. Then
\begin{itemize}
\item[\mbox{\rm(a)}] $bcirc(\mathscr{A}\ast\mathscr{B})=bcirc(\mathscr{A})bcirc(\mathscr{B})$,
\item[\mbox{\rm(b)}] $bcirc(\mathscr{A})^k=bcirc(\mathscr{A}^k)$ for all $k\in \{0,1,2,\ldots\}$,
\item[\mbox{\rm(c)}] $bcirc(\mathscr{A}^\top)=bcirc(\mathscr{A})^\top, bcirc(\mathscr{A}^H)=bcirc(\mathscr{A})^H$.
\end{itemize}
\end{lemma}

In the rest of this paper, for simplicity, we will use the following notation: for any $\mathscr{A}_i\in \mathbb{R}^{m_i\times n_i\times p}$ and $\mathscr{V}_i\in \mathbb{R}^{m_i\times n\times p}$ with $i\in [l]$, we denote
\begin{eqnarray*}\label{diag+vec}
\begin{array}{c}
Diag(\mathscr{A}_1,\mathscr{A}_2,\cdots,\mathscr{A}_l):= \left[\begin{array}{cccc}
\mathscr{A}_1&  &  &  \\
  & \mathscr{A}_2 &   &   \\
  &   & \ddots &  \\
  &   &   & \mathscr{A}_l
\end{array}\right],\;\;
vec(\mathscr{V}_1,\mathscr{V}_2,\cdots,\mathscr{V}_l):=\begin{bmatrix}
\mathscr{V}_{1}  \\
\mathscr{V}_{2}  \\
\vdots \\
\mathscr{V}_{l}
\end{bmatrix},
\end{array}
\end{eqnarray*}
and sometimes, they are abbreviated as $Diag(\mathscr{A}_i: i\in [l])$ and $vec(\mathscr{V}_i: i\in [l])$, respectively. When all $\mathscr{A}_i\; (\mathscr{V}_i)$ become matrices (or vectors or scalars), similar symbols are also used.

Recall that each circular matrix $\mathbf{A}\in \mathbb{R}^{n\times n}$ can be diagonalized with the normalized discrete Fourier transform (DFT) matrix \cite{GV-13}, i.e., there exists a diagonal matrix $\mathbf{D}$ such that $\mathbf{A}=\mathbf{F}_{n}^H\mathbf{D}\mathbf{F}_{n}$, where $\mathbf{F}_{n}$ is the Fourier matrix of size $n\times n$ defined as
\begin{eqnarray}\label{DFT}
\mathbf{F}_{n}=\frac{1}{\sqrt{n}}\left[
\begin{array}{cccccc}
  1 & 1 & 1 & 1 & \ldots & 1 \\
  1 & \omega & \omega^2 & \omega^3 & \cdots & \omega^{n-1} \\
  \vdots & \vdots & \vdots & \vdots & \ddots & \vdots \\
  1 & \omega^{n-1} & \omega^{2(n-1)} & \omega^{3(n-1)} & \cdots & \omega^{(n-1)(n-1)}
\end{array}
\right]
\end{eqnarray}
where $\omega=e^{\frac{2\pi \mathbf{i}}{n}}$ with $\mathbf{i}^2=-1$. Similarly, block circular matrices can be block diagonalized by using the Fourier transform. In \cite{KBHH-13}, the authors showed that for any third-order tensor $\mathscr{A}\in \mathbb{R}^{m\times n\times p}$, there exists a block diagonal matrix $Diag(\mathbf{A}_i: i\in [p])$ such that $bcirc(\mathscr{A})=(\mathbf{F}_{n}^H\otimes \mathbf{I}_{m\times m})Diag(\mathbf{A}_i: i\in [p])(\mathbf{F}_{p}\otimes\mathbf{I}_{n\times n})$, and pointed out the conjugate symmetry between these block matrices $\mathbf{A}_i$.

\begin{lemma}{\rm\cite{KBHH-13}}\label{coro-*}
Let $\mathscr{A}\in \mathbb{R}^{m\times n\times p}$ be block diagonalized as
\begin{eqnarray}\label{A0}
bcirc(\mathscr{A})=(\mathbf{F}_{p}^{H}\otimes\mathbf{I}_{m\times m})Diag(
\mathbf{A}_i: i\in [p])(\mathbf{F}_{p}\otimes\mathbf{I}_{n\times n}),
\end{eqnarray}
where $\mathbf{F}_{p}$ is the Fourier matrix defined by {\rm(\ref{DFT})}. Then, $\mathbf{A}_1\in \mathbb{R}^{m\times n}$, $\mathbf{A}_i\in \mathbb{C}^{m\times n}$ and $\mathbf{A}_i=\overline{\mathbf{A}_{p+2-i}}$ for any $i\in [p]\setminus\{1\}$. In particular, if $\mathscr{A}\in S\mathbb{R}^{n\times n\times p}$ and $\mathbf{A}^{(i)}\in S\mathbb{R}^{n\times n}$,
then $\mathbf{A}_i\in \mathbb{R}^{n\times n}$ for any $i\in [p]$ and $\mathbf{A}_i={\mathbf{A}_{p+2-i}}$ for any $i\in [p]\setminus\{1\}$.
\end{lemma}

\begin{remark}\label{remark-1}
It should be noticed that, for any $\mathscr{A}\in \mathbb{R}^{m\times n\times p}$ which can be block diagonalized as (\ref{A0}), most of the matrices $\mathbf{A}_i$ $(i\in [p])$ may be complex because of the participating of Fourier matrix $\mathbf{F}_p$, and they satisfy the relationships:
$\mathbf{A}_1\in \mathbb{R}^{m\times n}$ and $\mathbf{A}_i=\overline{\mathbf{A}_{p+2-i}}$ for any $i\in [p]\setminus\{1\}$. It should be noted that most of the matrices $\mathbf{A}_i$ $(i\in [p])$ may be complex even when $\mathscr{A}$ is symmetric, since $\mathbf{A}^{(i)}=\mathbf{A}^{(p+2-i)}$ for any $i\in [p]\setminus\{1\}$ may not hold in this case. On the other hand, it should be noticed that any $\mathbf{A}_i\in \mathbb{C}^{m\times n}$ with $i\in [p]$, which satisfy the above relationships, can lead to a real tensor by the construction as (\ref{A0}).
\end{remark}

\section{T-Hessian tensor and convexity of the multi-vector real-valued function}
\label{Sect. 2}
As is well-known to us, the local curvature of a multi-variable real-valued function can be characterized by the positive semidefiniteness of its Hessian matrix, which is widely used in Newton-type methods for solving various optimization problems.
In this section, we generalize the Hessian matrix to the third-order tensor: T-Hessian tensor, and study the discriminant condition of the convexity of the multi-vector real-valued function.

\subsection{Derivatives of multi-vector real-valued functions}\label{Sect. 3.1}
In this subsection, we explore the derivatives of multi-vector real-valued functions. The tensor space is a normed linear space with inner product. In the following, we regard a matrix as a third-order tensor and derive the multi-vector real-valued function with the help of the inner product and tensor T-product. We adopt the  Fr\'{e}chet derivative: Let $V$ and $W$ be normed vector spaces, and {$U\subseteq V$} be an open subset of $V$. A function $f : U\rightarrow W$ is called to be Fr\'{e}chet differentiable at $\mathbf{x}\in U$ if there exists a bounded linear operator $A:V\rightarrow W$ such that
$$\lim\limits_{\mathbf{h}\rightarrow \mathbf{0}}\frac{\|f(\mathbf{x}+\mathbf{h})- f(\mathbf{x})-A(\mathbf{h})\|_{W}}{ \|\mathbf{h}\|_{V} }=0.$$

Recently, in \cite{B-10} and \cite{KBHH-13}, the authors showed that third-order tensors can be seen as linear operators on a space of matrices with the help of the newly proposed tensor T-product and obtained many good theoretical and computational results. Inspired by that, we wonder whether or not can we adopt the third-order tensor as the linear operator in the above definition of Fr\'{e}chet derivative when the variety is a matrix? To this end, we introduce the following definition first.

\begin{definition}
\label{T-derivatives}
Let $f: U\subseteq\mathbb{R}^{n\times 1\times p}\rightarrow \mathbb{R}$ be a continuous map. Then, we say $f$ is T-differentiable at $\mathscr{X}\in U$ if and only if there exists a third-order tensor $\nabla_\mathscr{T} f(\mathscr{X})$ such that
$$\lim\limits_{\mathscr{H}\rightarrow \mathscr{O}}\frac{\|f(\mathscr{X}+\mathscr{H})- f(\mathscr{X})-\langle\nabla_\mathscr{T} f(\mathscr{X}),\mathscr{H}\rangle\|}{ \|\mathscr{H}\| }=0,$$
and we denote the first-order T-derivative of $f$ at $\mathscr{X}$ as $\nabla_\mathscr{T} f(\mathscr{X}):=\frac{\partial f}{\partial \mathscr{X}}$; and we say $f$ is twice T-differentiable at $\mathscr{X}\in U$ if and only if $f$ is continuously T-differentiable and there exists a third-order tensor $\nabla^2_\mathscr{T} f(\mathscr{X})$ such that
$$\lim\limits_{\mathscr{H}\rightarrow \mathscr{O}}\frac{\|\nabla_\mathscr{T} f(\mathscr{X}+\mathscr{H})- \nabla_\mathscr{T} f(\mathscr{X})-\nabla^2_\mathscr{T} f(\mathscr{X})\ast\mathscr{H}\|}{ \|\mathscr{H}\| }=0,$$
and we denote the second-order T-derivative of $f$ at $\mathscr{X}$ as $\nabla^2_\mathscr{T} f(\mathscr{X}):=\frac{\partial \nabla_\mathscr{T} f(\mathscr{X})}{\partial \mathscr{X}}$.

Furthermore, we say $f$ is T-differentiable (twice T-differentiable) on $U$ if and only if $f$ is T-differentiable (twice T-differentiable) at every $\mathscr{X}\in U$.
\end{definition}

\begin{remark}
From the fact that the tensor T-product of two tensors of size $m\times n\times p$ reduces to the matrix multiplication when $p=1$, it follows that $f: U\subseteq\mathbb{R}^{n\times 1\times p}\rightarrow \mathbb{R}$ being T-differentiable (twice T-differentiable) on $U$ is equivalent to $f$ being differentiable (twice differentiable) on $U$ when $p=1$.
\end{remark}

\begin{theorem}\label{general-derivative}
Let $f$ be a continuous map from $U\subseteq\mathbb{R}^{n\times 1\times p}$ to $\mathbb{R}$. Then
\begin{itemize}
\item[\mbox{\rm(i)}] $f$ is T-differentiable on $U$ if and only if $\frac{\partial f(\mathscr{X})}{\partial[unfold(\mathscr{X})]}$ exists for every $\mathscr{X}\in U$. Especially, for any $\mathscr{X}\in U$, $\nabla_\mathscr{T} f(\mathscr{X})=fold[\frac{\partial f(\mathscr{X})}{\partial[unfold(\mathscr{X})]}]$;\\
\item[\mbox{\rm(ii)}] $f$ is twice T-differentiable on $U$ if and only if $f$ is continuously T-differentiable on $U$ and $\frac{\partial [unfold(\nabla_\mathscr{T} f(\mathscr{X})]}{\partial [unfold (\mathscr{X})]}$ is a block circular matrix with each block being of size $n\times n$ for every $\mathscr{X}\in U$. In particular, $\nabla^2_\mathscr{T} f((\mathscr{X}))=bcirc^{-1}[\frac{\partial[ unfold(\nabla_\mathscr{T} f(\mathscr{X}))]}{\partial [unfold (\mathscr{X})]}]$ for any $\mathscr{X}\in U$.
\end{itemize}
\end{theorem}

\begin{proof} (i) $(\Rightarrow)$: If $f$ is T-differentiable on $U$, then for any $\mathscr{X}\in U$, there exists a bounded operator $\mathscr{L}$ such that
$$\lim\limits_{\mathscr{H\rightarrow\mathscr{O}}}\frac{\|f(\mathscr{X}+\mathscr{H})-f(\mathscr{X})-\langle\mathscr{L},\mathscr{H}\rangle\|}{\|\mathscr{H}\|}=0,$$
which, together with $\langle\mathscr{L},\mathscr{H}\rangle=\langle unfold(\mathscr{L}),unfold(\mathscr{H})\rangle$, implies that
$$\lim\limits_{unfold(\mathscr{H})\rightarrow\mathbf{O}}\frac{\|f(\mathscr{X}+\mathscr{H})-f(\mathscr{X})-\langle unfold(\mathscr{L}),unfold(\mathscr{H})\rangle\|}{\|unfold(\mathscr{H})\|}=0.$$
Furthermore, by introducing $g(unfold(\mathscr{X})):=f(\mathscr{X})$, we have that
$$\lim\limits_{\mathscr{H}\rightarrow\mathscr{O}}\frac{\|g[unfold(\mathscr{X}+\mathscr{H})]-g(unfold(\mathscr{X}))-\langle unfold(\mathscr{L}),unfold(\mathscr{H})\rangle\|}{\|unfold(\mathscr{H})\|}=0,$$
which means that $\frac{\partial [g(unfold(\mathscr{X}))]}{\partial[unfold({\mathscr{X}})]}=unfold(\mathscr{L})$, i.e., $\frac{\partial  [f(\mathscr{X})]}{\partial[unfold({\mathscr{X}})]}=unfold(\mathscr{L})$. Thus, $\frac{\partial f(\mathscr{X})}{\partial[unfold(\mathscr{X})]}$ exists, and $\nabla_\mathscr{T} f(\mathscr{X})=\mathscr{L}=fold[\frac{\partial f(\mathscr{X})}{\partial[unfold(\mathscr{X})]}]$.

$(\Leftarrow)$: Reversing the above process, we can obtain that if $\frac{\partial f(\mathscr{X})}{\partial[unfold(\mathscr{X})]}$ exists, then $f$ is T-differentiable on $U$.

(ii) $(\Rightarrow)$: If $f$ is twice T-differentiable on $U$,  then $f$ is continuously T-differentiable and for any $\mathscr{X}\in U$, there exists a bounded operator $\mathscr{L}: \mathbb{R}^{n\times 1\times p}\rightarrow\mathbb{R}^{n\times 1\times p}$ such that
$$\lim\limits_{\mathscr{H\rightarrow\mathscr{O}}}\frac{\|\nabla_\mathscr{T} f(\mathscr{X}+\mathscr{H})-\nabla_\mathscr{T} f(\mathscr{X})-\mathscr{L}\ast\mathscr{H}\|}{\|\mathscr{H}\|}=0,$$
which, together with $\|\mathscr{X}\|=\|unfold(\mathscr{X})\|$ for any $\mathscr{X}\in\mathbb{R}^{n\times 1\times p}$, implies that
$$\lim\limits_{\mathscr{H}\rightarrow\mathscr{O}}\frac{\|unfold[\nabla_\mathscr{T} f(\mathscr{X}+\mathscr{H})]-unfold[\nabla_\mathscr{T} f(\mathscr{X})]-bcirc(\mathscr{L})unfold(\mathscr{H})\|}{\|unfold(\mathscr{H})\|}=0.$$
Denote $h:\mathbb{R}^{np}\rightarrow \mathbb{R}^{np}$ with $h[unfold (\mathscr{X})]=unfold(\nabla_\mathscr{T} f(\mathscr{X}))$, then we have that $\frac{\partial [unfold(\nabla_\mathscr{T} f(\mathscr{X})]}{\partial [unfold (\mathscr{X})]}=\frac{\partial [h[unfold(\mathscr{X})]}{\partial [unfold (\mathscr{X})]}=bcirc(\mathscr{L})$, which is a block circular matrix with each block being of size $n\times n$. Thus $\frac{\partial [unfold(\nabla_\mathscr{T} f(\mathscr{X})]}{\partial [unfold (\mathscr{X})]}$ exists, and $$\nabla^2_\mathscr{T} f(\mathscr{X})=\mathscr{L}=bcirc^{-1}[\frac{\partial[ unfold(\nabla_\mathscr{T} f(\mathscr{X}))]}{\partial [unfold (\mathscr{X})]}].$$

$(\Leftarrow)$: Reversing the above process, we can obtain the desired result.
\end{proof}

In Theorem \ref{general-derivative}, we establish the necessary and sufficient conditions for a general multi-vector real-valued function being T-differentiable and twice T-differentiable, respectively, and we give the specific expressions when the derivatives exist. Then what are the relationships between the derivatives obtained in this way and the derivatives in the traditional sense? We construct an example to illustrate that.

{\bf Example 3.1.}
Given a map $f:\mathbb{R}^{2\times1\times 2}\rightarrow \mathbf{R}$ with for any $\mathscr{X}=[x_{i1j}]\in \mathbb{R}^{2\times1\times 2}$,
$$\begin{array}{ll}
f(\mathscr{X})&=x_{111}^2 + 2x_{111}x_{112} + x_{112}^2 + x_{211}^2
+ x_{212}^2. \end{array}$$

(1) First, we discuss the relationship between $\nabla f$ and $\nabla_\mathscr{T} f$.
By the traditional way, we can obtain that
$$\begin{array}{ll}
\nabla f(\mathscr{X})=
\begin{bmatrix}
  \frac{\partial f}{\partial x_{111}} &\frac{\partial f}{\partial x_{112}}\\
  \frac{\partial f}{\partial x_{211}} & \frac{\partial f}{\partial x_{212}}
\end{bmatrix}&
=\begin{bmatrix}
   2x_{111} + 2x_{112} & 2x_{111} + 2x_{112} \\
   2x_{211} & 2x_{212}
\end{bmatrix}
\end{array}.$$
Now, let us to seek the $\nabla_\mathscr{T} f$ by the procedure given in Theorem \ref{general-derivative}. Noting that
$unfold({\mathscr{X}})=
[{x_{111}}, {x_{211}} ,  {x_{112}},{ x_{212}}]^\top,$
thus we can obtain that
$$\begin{array}{l}
\frac{\partial  [f(\mathscr{X})]}{\partial[unfold({\mathscr{X}})]}=
[\frac{\partial f}{\partial x_{111}} ,\frac{\partial f}{\partial x_{211}} ,\frac{\partial f}{\partial x_{112}}, \frac{\partial f}{\partial x_{212}}]^\top=2[
x_{111} + x_{112}, x_{211}, x_{111} + x_{112}, x_{212}]^\top,
\end{array}$$
furthermore, we can get
$$\begin{array}{l}
\nabla_\mathscr{T} f(\mathscr{X})=fold[\frac{\partial f(\mathscr{X})}{\partial[unfold(\mathscr{X})]}]=\begin{bmatrix}
   2x_{111} + 2x_{112} & 2x_{111} + 2x_{112} \\
   2x_{211} & 2x_{212}
\end{bmatrix}=\nabla f(\mathscr{X}).
\end{array}$$

(2) Next, we investigate the relationship between $\nabla^2 f$ and $\nabla^2_\mathscr{T} f$.

Since $\nabla f(\mathscr{X})$ is actually a matrix of size $2\times 2$ and $\mathscr{X}$ is also a matrix of size $2\times 2$, thus by the traditional derivative, $\nabla^2 f(\mathscr{X})$ should be a tensor of fourth-order with the elements being same as ones in the matrix $\mathbf{A}\in \mathbb{R}^{4\times 4}$ where:
$$\mathbf{A}=\frac{\partial (\nabla f(\mathscr{X}))}{\partial \mathscr{X}}=[\frac{\partial (\nabla f(\mathscr{X}))}{\partial x_{ij}}]_{4\times 4}=
\begin{bmatrix}
\begin{array}{cccc}
  2 & 2 & 2 & 2 \\
  0 & 0 & 0 & 0 \\
  0 & 0 & 0 & 0 \\
  2 & 0 & 0 & 2
\end{array}
\end{bmatrix}.$$

While, by Theorem \ref{general-derivative}, we can get that
$$bcirc(\nabla^2_\mathscr{T} f(\mathscr{X}))=\frac{\partial[ unfold(\nabla_\mathscr{T} f(\mathscr{X}))]}{\partial [unfold (\mathscr{X})]}=
\begin{bmatrix}
\begin{array}{cccc}
  2 & 0 & 2 & 0 \\
  0 & 2 & 0 & 0 \\
  2 & 0 & 2 & 0 \\
  0 & 0 & 0 & 2
\end{array}
\end{bmatrix}$$
which is a block circular matrix and $\nabla^2_\mathscr{T} f(\mathscr{X})$ is a tensor of size $2\times 2\times 2$, with the frontal slices being:
$$\begin{array}{ccc}
(\nabla^2_\mathscr{T}f(\mathscr{X}))^{(1)}=\begin{bmatrix}
\begin{array}{cc}
  2 & 0 \\
  0 & 2
\end{array}
\end{bmatrix}& \mbox{\rm and} &(\nabla^2_\mathscr{T}f(\mathscr{X}))^{(2)}=\begin{bmatrix}
\begin{array}{cc}
  2 & 0 \\
  0 & 0
\end{array}
\end{bmatrix}
\end{array}.$$

Hence, $\nabla^2 f$ and $\nabla^2_{\mathscr{T}} f$ are not coincide in the sizes, but it should be noticed that the entries of $\nabla^2 f$ and these of $bcirc(\nabla^2_{\mathscr{T}}) f$ are the same regardless of the orders.

{\begin{remark}
(i) $\nabla^2_\mathscr{T} f(\mathscr{X})$ is different from the traditional one due to the participation of the T-product. Traditionally, the second-order derivative of a multi-vector real-valued function is a fourth-order tensor. However, $\nabla^2_\mathscr{T} f(\mathscr{X})$ turns out to be a third-order tensor, which is also reasonable. Its rationality lies in that the existence of $\nabla^2_\mathscr{T} f(\mathscr{X})$ via the T-product needs such a condition that $\frac{\partial [unfold(\nabla f(\mathscr{X})]}{\partial [unfold (\mathscr{X})]}$ is a block circular matrix, which implies that just getting the information of its first block column vector is enough. In other words, it is enough to express the information of $\nabla^2_\mathscr{T} f(\mathscr{X})$ by a third-order tensor in such case.


(ii) The $\nabla_\mathscr{T} f(\mathscr{X})$ is consistent with the traditional one.
This is natural because $\nabla_\mathscr{T} f(\mathscr{X})$ and $\nabla f(\mathscr{X})$ are based on the coincide definition of inner product. So, we use $\nabla f(\mathscr{X})$ to represent $\nabla_\mathscr{T} f(\mathscr{X})$ in the rest of paper.

\end{remark}}

\subsection{ The semidefiniteness of $\nabla^2_\mathscr{T} f(\mathscr{X})$ and the convexity of $f(\mathscr{X})$}\label{Sect. 3.2}
In this subsection, we investigate the second-order condition for any twice continuously T-differentiable function $f:U\subseteq \mathbb{R}^{n\times1\times p}\rightarrow\mathbb{R}$ being (strictly) convex. The definition of the convex matrix function is as follows.
\begin{definition}\cite{B-04}
A function $f: U\subseteq\mathbb{R}^{n\times p}\rightarrow\mathbb{R}$ is convex (strictly convex) if $U$ is a convex set, and for all $\mathbf{X},\mathbf{Y}\in U$ ($\mathbf{X},\mathbf{Y}\in U$ and $\mathbf{X}\neq\mathbf{Y}$) and any $\theta$ with $0\leq\theta\leq 1$,
$
f(\theta \mathbf{X}+(1-\theta) \mathbf{Y})\leq(<) \theta f(\mathbf{X})+(1-\theta)f (\mathbf{Y}).
$
\end{definition}

Since the first-order T-derivative is consistent with the traditional one for a multi-vector real-valued function, it is not difficult to find that for any continuously T-differentiable function $f:U\subseteq \mathbb{R}^{n\times1\times p}\rightarrow\mathbb{R}$, it is convex (strictly convex) if and only if for any $\mathscr{X},\mathscr{Y}\in U$ ($\mathscr{X},\mathscr{Y}\in U$ and $\mathscr{X}\neq\mathscr{Y}$),
$$f(\mathscr{Y})\geq(>) f(\mathscr{X})+ \langle\nabla f(\mathscr{X}),\mathscr{Y}-\mathscr{X}\rangle.$$
Below, in order to establish the second-order condition for any twice continuously T-differentiable function $f:U\subseteq \mathbb{R}^{n\times1\times p}\rightarrow\mathbb{R}$ being (strictly) convex, we first extend the second-order Taylor expansion of the function of one variable to the multi-vector real-valued function via the tensor T-product.

\begin{theorem}\label{Taylor expansion}
Suppose that $f:U\subseteq \mathbb{R}^{n\times1\times p}\rightarrow\mathbb{R}$ is twice continuously T-differentiable on $U$. Then
\begin{itemize}
\item[\mbox{\rm(i)}] there exists $\beta\in(0,1)$ such that
$$f(\mathscr{X})= f(\widetilde{\mathscr{X}})+ \langle\nabla f(\widetilde{\mathscr{X}}),\mathscr{X}-\widetilde{\mathscr{X}}\rangle+\frac{1}{2}\langle\nabla^2_\mathscr{T} f(\mathscr{Z})\ast(\mathscr{X}-\widetilde{\mathscr{X}}),\mathscr{X}-\widetilde{\mathscr{X}}\rangle,
$$
where $\mathscr{Z}=\beta \mathscr{X}+ (1-\beta)\widetilde{\mathscr{X}}$;
\item[\mbox{\rm(ii)}] the second-order T-Taylor expansion of $f$ at $\widetilde{\mathscr{X}}$ as follows:
\begin{eqnarray*}
\begin{array}{ll}
f(\mathscr{X})=& f(\widetilde{\mathscr{X}})+ \langle\nabla f(\widetilde{\mathscr{X}}),\mathscr{X}-\widetilde{\mathscr{X}}\rangle\\
&+\frac{1}{2}\langle\nabla^2_\mathscr{T} f(\widetilde{\mathscr{X}})\ast(\mathscr{X}-\widetilde{\mathscr{X}}),\mathscr{X}-\widetilde{\mathscr{X}}\rangle+o(\|\mathscr{X}-\widetilde{\mathscr{X}}\|^2)
\end{array}\end{eqnarray*}
and $o(\|\mathscr{X}-\widetilde{\mathscr{X}}\|^2)$ means a high-order infinitesimal of $\|\mathscr{X}-\widetilde{\mathscr{X}}\|^2$ as $\mathscr{X}\rightarrow \widetilde{\mathscr{X}}$.
\end{itemize}
\end{theorem}

\begin{proof} (i) Construct a function of one variable as: $\varphi(t):=f(\widetilde{\mathscr{X}}+t(\mathscr{X}-\widetilde{\mathscr{X}}))$ for any $t\in \mathbb{R}$. Then, $\varphi(0)=f(\widetilde{\mathscr{X}})$
and $\varphi(1)=f({\mathscr{X}})$. It follows from the condition that $f$ is twice continuously T-differentiable that $\varphi$ is twice continuously differentiable. For any $t\in \mathbb{R}$, let $\mathscr{U}=\widetilde{\mathscr{X}}+t(\mathscr{X}-\widetilde{\mathscr{X}})$, then we have
$$\renewcommand\arraystretch{1.2}
\begin{array}{lcl}
\varphi'(t)&=&\frac{d\varphi(t)}{dt}=\frac{d f(\mathscr{U})}{dt}=\langle\nabla f(\mathscr{U}),\frac{d\mathscr{U}}{dt}\rangle=\langle\nabla f(\mathscr{U}),(\mathscr{X}-\widetilde{\mathscr{X}})\rangle\\
&=& \langle\nabla f(\widetilde{\mathscr{X}}+t(\mathscr{X}-\widetilde{\mathscr{X}})),(\mathscr{X}-\widetilde{\mathscr{X}})\rangle;\\
\varphi''(t)&=& \frac{d\varphi'(t)}{dt}=\frac{d\langle\nabla f(\mathscr{U}),(\mathscr{X}-\widetilde{\mathscr{X}})\rangle}{dt}=\langle\frac{d(\nabla f(\mathscr{U}))}{dt},(\mathscr{X}-\widetilde{\mathscr{X}})\rangle\\
&=&\langle\frac{\nabla^2_\mathscr{T} f(\mathscr{U})\ast d\mathscr{U}}{dt},\mathscr{X}-\widetilde{\mathscr{X}}\rangle\\
&=&\langle\frac{bcirc(\nabla^2_\mathscr{T}f(\mathscr{U}))\cdot unfold(d\mathscr{U})}{dt},unfold(\mathscr{X}-\widetilde{\mathscr{X}})\rangle\\
&=&\langle bcirc(\nabla^2_\mathscr{T}f(\mathscr{U}))\cdot unfold(\mathscr{X}-\widetilde{\mathscr{X}}),unfold(\mathscr{X}-\widetilde{\mathscr{X}})\rangle\\
&=&\langle unfold(\nabla^2_\mathscr{T}  f(\mathscr{U})\ast(\mathscr{X}-\widetilde{\mathscr{X}})),unfold(\mathscr{X}-\widetilde{\mathscr{X}})\rangle\\
&=&\langle\nabla^2_\mathscr{T}  f(\widetilde{\mathscr{X}}+t(\mathscr{X}-\widetilde{\mathscr{X}}))\ast(\mathscr{X}-\widetilde{\mathscr{X}}), \mathscr{X}-\widetilde{\mathscr{X}}\rangle.
\end{array}$$
Thus, $\varphi'(0)=\langle\nabla  f(\widetilde{\mathscr{X}}),\mathscr{X}-\widetilde{\mathscr{X}}\rangle$ and $\varphi''(0)=\langle\nabla^2_\mathscr{T}  f(\widetilde{\mathscr{X}})\ast(\mathscr{X}-\widetilde{\mathscr{X}}), \mathscr{X}-\widetilde{\mathscr{X}}\rangle$.
It follows from the mean value theorem, that there exists some $\beta\in (0,1)$ such that $\varphi(1)=\varphi(0)+\varphi'(0)+\frac{1}{2}\varphi''(\beta)$, which implies that
$$
\begin{array}{rcl}
f(\mathscr{X})&=& f(\widetilde{\mathscr{X}})+ \langle\nabla f(\widetilde{\mathscr{X}}),\mathscr{X}-\widetilde{\mathscr{X}}\rangle+\frac{1}{2}\langle\nabla^2_\mathscr{T} f(\mathscr{Z})\ast(\mathscr{X}-\widetilde{\mathscr{X}}),\mathscr{X}-\widetilde{\mathscr{X}}\rangle,
\end{array}$$
where
$\mathscr{Z}=\beta \mathscr{X}+ (1-\beta)\widetilde{\mathscr{X}}$, i.e., the result in (i) holds.

(ii) Denote $\mathscr{Y}:=\frac{\mathscr{X}-\widetilde{\mathscr{X}}}{\|\mathscr{X}-\widetilde{\mathscr{X}}\|}$ and $\gamma':=\|\mathscr{X}-\widetilde{\mathscr{X}}\|$. Let $\psi(\gamma):=f(\widetilde{\mathscr{X}}+\gamma\mathscr{Y})$ for any $\gamma\in \mathbb{R}$, then by the same way as $(i)$, we have $\psi(0)=f(\widetilde{\mathscr{X}})$, $\psi(\gamma')=f({\mathscr{X}})$,
$$
\psi'(0)\gamma'=\langle\nabla  f(\widetilde{\mathscr{X}}),\mathscr{X}-\widetilde{\mathscr{X}}\rangle\;\;
\mbox{\rm and}\;\; \psi''(0)\gamma'^2=\langle\nabla^2_\mathscr{T}  f(\widetilde{\mathscr{X}})\ast(\mathscr{X}-\widetilde{\mathscr{X}}), \mathscr{X}-\widetilde{\mathscr{X}}\rangle,
$$
which, together with $\psi(\gamma')=\psi(0)+\psi'(0)\gamma'+\frac{1}{2}\psi''(0)\gamma'^2+o(\|\gamma'\|^2)$, imply that
$$\begin{array}{lll}
f(\mathscr{X})&=& f(\widetilde{\mathscr{X}})+ \langle\nabla f(\widetilde{\mathscr{X}})\ast(\mathscr{X}-\widetilde{\mathscr{X}}),\mathscr{X}-\widetilde{\mathscr{X}}\rangle\\
&&+\frac{1}{2}\langle\nabla^2_\mathscr{T} f(\widetilde{\mathscr{X}})\ast(\mathscr{X}-\widetilde{\mathscr{X}}),\mathscr{X}-\widetilde{\mathscr{X}}\rangle +o(\|\mathscr{X}-\widetilde{\mathscr{X}}\|^2),
\end{array}$$
i.e., the result in (ii) holds.
\end{proof}

\begin{theorem}\label{convex-condition}
Let $f: U\subseteq\mathbb{R}^{n\times1\times p}\rightarrow \mathbb{R}$ be a twice continuously T-differentiable function on an open convex set $U$. Then
\begin{itemize}
\item[\mbox{\rm(i)}] $f$ is convex if and only if for any $\mathscr{X}\in U$, $\nabla^2_\mathscr{T}f(\mathscr{X})$ satisfies
    \begin{eqnarray}\label{second-condition}
    \langle\mathscr{Y},\nabla^2_\mathscr{T}f(\mathscr{X})\ast \mathscr{Y}\rangle\geq 0 \;\; \mbox{\rm for\ any\ }\mathscr{Y}\in \mathbb{R}^{n\times1\times p};
    \end{eqnarray}
\item[\mbox{\rm(ii)}]
$f$ is strictly convex if for any $\mathscr{X}\in U$, $\nabla^2_\mathscr{T}f(\mathscr{X})$ satisfies
$$
\langle\mathscr{Y},\nabla^2_\mathscr{T}f(\mathscr{X})\ast \mathscr{Y}\rangle>0 \;\; \mbox{\rm for\ any\ }\mathscr{Y}\in \mathbb{R}^{n\times1\times p}\setminus\{\mathscr{O}\}.
$$
\end{itemize}
\end{theorem}

\begin{proof} (i) $(\Rightarrow)$: For any $\mathscr{X}\in U$ and $\mathscr{Y}\in \mathbb{R}^{n\times1\times p}\setminus \{\mathscr{O}\}$, it follows from $U$ being an open set that there exists $\varepsilon>0$ such that $\mathscr{X}+\theta\mathscr{Y}\in U$ when $\theta\in (-\varepsilon,\varepsilon)$. Since $f$ is convex, we have that $f(\mathscr{X}+\theta\mathscr{Y})\geq f(\mathscr{X})+\theta\langle\nabla f(\mathscr{X}),\mathscr{Y}\rangle$. In addition, it follows from Theorem \ref{Taylor expansion} that $$f(\mathscr{X}+\theta\mathscr{Y})=f(\mathscr{X})+\theta\langle\nabla f(\mathscr{X}),\mathscr{Y}\rangle+\frac{1}{2}\theta^2\langle\mathscr{Y},\nabla^2_\mathscr{T}f(\mathscr{X})\ast \mathscr{Y}\rangle+o(\|\theta\mathscr{Y}\|^2).$$ Therefore, we can obtain that
$\langle\mathscr{Y},\nabla^2_\mathscr{T}f(\mathscr{X})\ast \mathscr{Y}\rangle+\frac{o(\|\theta\mathscr{Y}\|^2)}{\theta^2/2}\geq 0$.
Let $\theta\rightarrow 0$, we further obtain that $\langle\mathscr{Y},\nabla^2_\mathscr{T}f(\mathscr{X})\ast \mathscr{Y}\rangle\geq 0$.

$(\Leftarrow)$: For any $\mathscr{X},\mathscr{Y}\in U$, it follows from Theorem \ref{Taylor expansion} that
$$f(\mathscr{Y})=f(\mathscr{X})+\langle\nabla f(\mathscr{X}),\mathscr{Y}-\mathscr{X}\rangle+\frac{1}{2}\langle\nabla^2_\mathscr{T}f(\mathscr{Z})\ast (\mathscr{Y}-\mathscr{X}),\mathscr{Y}-\mathscr{X}\rangle,$$
where $\mathscr{Z}=\mathscr{X}+t(\mathscr{Y}-\mathscr{X})$ with $t\in (0,1)$. Since $U$ is convex, it follows that $\mathscr{Z}\in U$; and hence, by (\ref{second-condition}) we have that $\frac{1}{2}\langle\mathscr{Y},\nabla^2_\mathscr{T}f(\mathscr{Z})\ast \mathscr{Y}\rangle\geq 0$. Furthermore, we have that $f(\mathscr{Y})\geq f(\mathscr{X})+\langle\nabla f(\mathscr{X}),\mathscr{Y}-\mathscr{X}\rangle$, which, together with $U$ being convex, implies that $f$ is convex. The proof of (i) is complete.

By the same way as in the proof of the sufficiency of (i), we can obtain (ii). \end{proof}

\begin{remark}
Let $f: U\subseteq\mathbb{R}^{n\times 1\times p}\rightarrow \mathbb{R}$ be a twice continuously T-differentiable function on an open convex set $U$. Since $\nabla^2_\mathscr{T}f\in \mathbb{R}^{n\times n\times p}$ has similar properties as Hessian matrix, we call $\nabla^2_\mathscr{T}f$ the T-Hessian tensor.
\end{remark}

\section{Symmetric T-positive (semi)definite tensors}\label{Sect. 3}

In this section, we first introduce a definition of the symmetric T-positive (semi)definite tensor; and then, we investigate properties of symmetric T-positive (semi)definite tensors.

\subsection{ Definition}\label{Sect. 4.1} In Section \ref{Sect. 2}, we obtained that the convexity of a twice continuously T-differentiable multi-vector real-valued function on an open convex set can be characterized by some property of the T-Hessian tensor. Now we name such a property as the symmetric T-positive semidefiniteness.

\begin{definition}\label{STPSD}
Let $\mathscr{A}\in \mathbb{R}^{n\times n\times p}$. We say $\mathscr{A}$ is a symmetric T-positive (semi)definite tensor {(T-P(S)D tensor for short)}, if and only if $\mathscr{A}$ is a symmetric tensor and $$\langle\mathscr{X}, \mathscr{A}\ast \mathscr{X}\rangle > (\geq) 0$$
holds for any $\mathscr{X}\in \mathbb{R}^{n\times 1\times p}\setminus\{\mathscr{O}\}$ (for any $\mathscr{X}\in \mathbb{R}^{n\times 1\times p}$). We denote the set consisting of all symmetric T-P(S)D tensors of size $n\times n\times p$ as $S\mathbb{R}^{n\times n\times p}_{++}(S\mathbb{R}^{n\times n\times p}_+)$.
\end{definition}

\begin{remark}
(i) When $p=1$, the T-product defined by Definition \ref{T-product} reduces to the product of two matrices. In addition, when $p=1$, $\mathscr{X}\in \mathbb{R}^{n\times 1\times p}$ reduces to a column vector and $\mathscr{A}\in S\mathbb{R}^{n\times n\times p}$ reduces to a square symmetric matrix. Thus, when $p=1$, Definition \ref{STPSD} is exactly the definition of the symmetric P(S)D matrix. That is to say, the T-P(S)D tensor defined by Definition \ref{STPSD} is a higher-order extension of the P(S)D matrix.

(ii) From Definition \ref{STPSD} and Theorem \ref{convex-condition}, we can say that a twice continuously T-differentiable function $f: U\subseteq\mathbb{R}^{n\times 1\times p}\rightarrow \mathbb{R}$ is convex (strictly convex) if and only if (if) the T-Hessian tensor $\nabla^2_\mathscr{T}f(\mathscr{X})$ is symmetric T-PSD (T-PD) for any $\mathscr{X}\in U$.

(iii) {It should be noticed that the positive semidefinite tensor defined by means of the nonnegativity of the corresponding multi-variate homogeneous polynomial in \cite{QL-17} is different with the one defined in Definition \ref{STPSD}. Since the positive semidefinite tensor defined by Qi \cite{QL-17} vanishes when the order is odd, while the symmetric T-positive (semi)definite tensor in Definition \ref{STPSD} is introduced for third-order tensor.}
\end{remark}

\subsection{ Equivalent characterizations of symmetric T-P(S)D tensors}\label{Sect 4.2}
First, we give an equivalent description of Definition \ref{STPSD}.
\begin{theorem}\label{equivalent}
Suppose that $\mathscr{A}\in \mathbb{R}^{n\times n\times p}$ can be block diagonalized as
\begin{eqnarray}\label{PSD-A}
bcirc(\mathscr{A})=(\mathbf{F}_{p}^{H}\otimes\mathbf{I}_{n\times n})Diag(
\mathbf{A}_i: i\in [p])
(\mathbf{F}_{p}\otimes\mathbf{I}_{n\times n}),
\end{eqnarray}
where $\mathbf{F}_{p}$ is the Fourier matrix of size $p\times p$, which is defined as (\ref{DFT}). Then $\mathscr{A}$ is symmetric T-P(S)D if and only if all the matrices $\mathbf{A}_i$ are Hermitian P(S)D.
\end{theorem}

\begin{proof} By (\ref{PSD-A}), $bcirc(\mathscr{A})$ is symmetric if and only if each $\mathbf{A}_i$ is Hermitian, and $\mathscr{A}$ is symmetric if and only if $bcirc(\mathscr{A})$ is symmetric as shown in \ref{pro-1(a)}.

$(\Leftarrow)$: Suppose that all the matrices $\mathbf{A}_i$ in (\ref{PSD-A}) are Hermitian P(S)D, then for any $\mathbf{x}$ in $\mathbb{C}^n\setminus\{\mathbf{0}\}$ ($\mathbf{x}\in \mathbb{C}^n$) and $i\in[p]$, $\mathbf{x}^H\mathbf{A}_i\mathbf{x}>(\geq)0$. For any $\mathscr{X}\in\mathbb{R}^{n\times 1 \times p}\setminus\{\mathscr{O}\}$, there exists $\mathbf{x}_{i}\in \mathbb{C}^{n}$ for each $i\in[p]$, which cannot be $\mathbf{0}$ at the same time, such that
$bcirc(\mathscr{X})=(\mathbf{F}_{p}^{H}\otimes\mathbf{I}_{n\times n})
Diag(\mathbf{x}_i: i\in [p])\mathbf{F}_{p}$.
Since
$$
\begin{array}{rcl}
\langle\mathscr{X}, \mathscr{A}\ast \mathscr{X}\rangle&=& \frac{1}{p}\langle bcirc(\mathscr{X}), bcirc(\mathscr{A})bcirc(\mathscr{X})\rangle\\
&=&\frac{1}{p}Tr(bcirc(\mathscr{X})^H bcirc(\mathscr{A})bcirc(\mathscr{X}))\\
&=&\frac{1}{p}Tr(\mathbf{F}_{p}^{H}Diag(\mathbf{x}^H_i\mathbf{A}_i\mathbf{x}_i: i\in [p])\mathbf{F}_{p})\\
&=&\frac{1}{p}Tr(Diag(\mathbf{x}^H_i\mathbf{A}_i\mathbf{x}_i: i\in [p]))
=\frac{1}{p}\sum\limits_{i=1}^p(\mathbf{x}^H_i\mathbf{A}_i\mathbf{x}_i)\geq 0,
\end{array}
$$
where the fourth equality follows from the fact that similar matrices have the same traces, then $\mathscr{A}$ is symmetric T-PSD if all the matrices $\mathbf{A}_i$ are Hermitian PSD.

Since $\mathbf{x}_{i}\in \mathbb{C}^{n}$ cannot be $\mathbf{0}$ at the same time, there exists at least an index $i\in [p]$ such that $\mathbf{x}^H_i\mathbf{A}_i\mathbf{x}_i>0$ if all the matrices $\mathbf{A}_i$ are Hermitian PD. Hence $\mathscr{A}$ is symmetric T-PD if all the matrices $\mathbf{A}_i$ are Hermitian PD.

$(\Rightarrow)$: Suppose that $\mathscr{A}\in S\mathbb{R}^{n\times n\times p}_{++}$ ($S\mathbb{R}^{n\times n\times p}_+$), then $\langle\mathscr{X}, \mathscr{A}\ast \mathscr{X}\rangle > (\geq) 0$ for any $\mathscr{X}\in\mathbb{R}^{n\times 1 \times p}\setminus\{\mathscr{O}\}$. Below, we divide the discussion into two cases.

{\bf Case 1: $n$ is even.} By Lemma \ref{coro-*}, we have $\mathbf{A}_1\in \mathbb{R}^{n\times n}$, $\mathbf{A}_{\frac{p+2}{2}}\in \mathbb{R}^{n\times n}$, $\mathbf{A}_i\in \mathbb{C}^{n\times n}$ and $\mathbf{A}_i=\overline{\mathbf{A}_{p+2-i}}$ for any $i\in
[p]\setminus\{1,\frac{p+2}{2}\}$. Then, for any $\mathbf{x}\in \mathbb{C}^{n}\setminus\{\mathbf{0}\}$, choose special $\mathscr{X'}_{i}$ in $\mathbb{R}^{n\times 1 \times p}$ and satisfies that
$
bcirc(\mathscr{X}')=(\mathbf{F}^{H}_{p}\otimes\mathbf{I}_{n\times n})Diag(\mathbf{X'}_i: i\in [p])\mathbf{F}_{p},
$
where $\mathbf{X'}_k=\overline{\mathbf{X}_{p+2-k}'}=\mathbf{x}$ and others $\mathbf{X'}_i=\mathbf{0}$ with $k$ being any fixed number in $[p]\setminus\{1,\frac{p+2}{2}\}$. Then, from Remark \ref{remark-1} we have that $\mathscr{X'}\in\mathbb{R}^{n\times 1\times p}$. Thus, $\langle\mathscr{X'},\mathscr{A}\ast\mathscr{X'})>(\geq 0)$ by $\mathscr{A}\in S\mathbb{R}^{n\times n\times p}_{++}(\mathscr{A}\in S\mathbb{R}^{n\times n\times p}_+)$. Since $\mathscr{A}\in S\mathbb{R}^{n\times n\times p}$, $\mathbf{A}_k\in H\mathbb{C}^{n\times n}$ and $\mathbf{A}_k=\overline{\mathbf{A}_{p+2-k}}$. Thus, $\mathbf{x}^H\mathbf{A}_k\mathbf{x}$ is real, and
$$0<(\leq) \langle\mathscr{X'},\mathscr{A}\ast\mathscr{X'})=\frac{1}{p}\sum\limits_{i=1}^p(\mathbf{x}^H_i\mathbf{A}_i\mathbf{x}_i)=
\frac{1}{p}(\mathbf{x}^H\mathbf{A}_k\mathbf{x}+ \overline{\mathbf{x}}^H\mathbf{A}_{n+2-k}\overline{\mathbf{x}})=\frac{2}{p} \mathbf{x}^H\mathbf{A}_k\mathbf{x},$$
which implies $\mathbf{A}_k$ is Hermitian P(S)D for any $k\in [p]\setminus\{1,\frac{p+2}{2}\}$. In addition, for any $\mathbf{x}\in \mathbb{R}^{n}$, we can obtain that $\mathbf{A}_1$ ($\mathbf{A}_{\frac{p+2}{2}}$) is symmetric P(S)D by choosing $\mathbf{X'}_1=\mathbf{x}$ ($\mathbf{X'}_{\frac{p+2}{2}}=\mathbf{x}$) and others $\mathbf{X'}_i=\mathbf{0}$.

{\bf Case 2: $n$ is odd.} By Lemma \ref{coro-*}, we have $\mathbf{A}_1\in \mathbb{R}^{n\times n}$, $\mathbf{A}_i\in \mathbb{C}^{n\times n}$ and $\mathbf{A}_i=\overline{\mathbf{A}_{p+2-i}}$ for any $i\in [p]\setminus\{1\}$. Then by the same method as Case 1, we can obtain that $\mathbf{A}_k$ is Hermitian P(S)D for any $k\in [p]$.

Thereby, combining {\bf{Case 1}} and {\bf{Case 2}}, we can obtain that all the matrices $\mathbf{A}_i$ are Hermitian P(S)D if $\mathscr{A}$ is symmetric T-P(S)D. \end{proof}

\begin{remark}
Theorem \ref{equivalent} shows that the judgement of the T-positive semidefiniteness of a symmetric tensor of size $n\times n\times p$ can be transformed into the judgement of positive semidefiniteness of $p$ Hermitian matrices of size $n\times n$. Furthermore, Theorem \ref{equivalent} shows that the symmetric T-P(S)D tensor in Definition \ref{STPSD} is equivalent to the one by \cite[Definition 2.7]{KBHH-13} and the one by \cite[Definition 15]{MQW-19} in real case.
\end{remark}

Next, we give another equivalent description of Definition \ref{STPSD}.
\begin{theorem}\label{pro-1(a)}
Let $\mathscr{A}\in S\mathbb{R}^{n\times n\times p}$. $\mathscr{A}$ is symmetric T-P(S)D if and only if $bcirc(\mathscr{A})$ is symmetric P(S)D.
\end{theorem}

\begin{proof} Since $bcirc(\mathscr{A}^\top)=bcirc(\mathscr{A}^H)=bcirc(\mathscr{A})^H=bcirc(\mathscr{A})^\top$ by Lemma \ref{lemma-1}$(c)$, then $\mathscr{A}$ is symmetric if and only if $bcirc(\mathscr{A})$ is symmetric. For any $\mathscr{X}\in \mathbb{R}^{n\times 1\times p}$, it follows from Definition \ref{T-product} and the definition of operator $unfold$ that $$unfold(\mathscr{A}\ast \mathscr{X})=bcirc(\mathscr{A})unfold(\mathscr{X}),$$ and hence,
$$\begin{array}{lcl}
\langle\mathscr{X}, \mathscr{A}\ast \mathscr{X}\rangle&=&\langle unfold(\mathscr{X}), unfold(\mathscr{A}\ast \mathscr{X})\rangle\\
&=&\langle unfold(\mathscr{X}),bcirc(\mathscr{A})unfold(\mathscr{X})\rangle.
\end{array}$$
Thus, by combining Definition \ref{STPSD} and the criterion of P(S)D matrix, we can easily obtain that $\mathscr{A}$ is symmetric T-P(S)D if and only if $bcirc(\mathscr{A})$ is symmetric P(S)D. \end{proof}

\begin{remark}
From Theorem \ref{pro-1(a)}, we can see that lots of results that P(S)D matrix with block circular structure hold is true for T-P(S)D tensors by combining the properties of tensor T-product such as those shown in Lemma \ref{lemma-1}. Thus for convenience, in Section \ref{Sect. 4.3} and Section \ref{Sect. 4.4} we just list some ones which play important roles in Section \ref{Sect. 4.5} and Section \ref{Sect. 4} without proofs.
\end{remark}

\subsection{T-eigenvalue decomposition of the symmetric T-P(S)D tensor}\label{Sect. 4.3}
In this subsection, we aim to establish the T-eigenvalue decomposition for the symmetric T-P(S)D tensor. To do this, we give the following definition first. It should be noticed that the definition of T-eigenvalue for third-order $F$-square tensor was given in \cite{MQW-19} and here we redefine it in an equivalent way for convenience.

\begin{definition}(T-eigenvalue and Trace)\label{Tr}
Let $\mathscr{A}\in \mathbb{R}^{n\times n\times p}$, which can be block diagonalized as (\ref{PSD-A}). Then a real number $\lambda$ is said to be a T-eigenvalue of $\mathscr{A}$ if and only if it is an eigenvalue of some $\mathbf{A}_i$ for $i\in [p]$, denoted by $\lambda(\mathscr{A})$. The largest and smallest T-eigenvalues of $\mathscr{A}$ are denoted by $\lambda_{\max}(\mathscr{A})$ and $\lambda_{\min}(\mathscr{A})$, respectively. Moreover, the trace of $\mathscr{A}$, denoted by $Tr(\mathscr{A})$, is defined as $Tr(\mathscr{A}):=\sum_{i=1}^pTr(\mathbf{A}_i)$.
\end{definition}

\begin{remark}\label{remark-3}
By Definition \ref{Tr}, Theorem \ref{equivalent} and \cite[Fact 6]{T-01}, it is not difficult to see that a symmetric third-order tensor $\mathscr{A}$ is T-P(S)D if and only if each T-eigenvalue of $\mathscr{A}$ is positive (nonnegative).
\end{remark}

It is easy to establish the following properties for the T-eigenvalues and traces of tensors from the above definition and some known results in \cite{HJ-13}.

\begin{proposition}\label{pro-2}
Let $\mathscr{A}$ and $\mathscr{B}$ be two tensors in $\mathbb{R}^{n\times n\times p}$, $\mathscr{C}\in \mathbb{R}^{n\times n\times p}$ be nonsingular, and $spec(\mathscr{A})$ be the set consisting of all the T-eigenvalues of $\mathscr{A}$. Then
\begin{itemize}
\item[(a)] $spec(\mathscr{A})=spec(bcirc(\mathscr{A}))$;
\item[(b)] $Tr(\mathscr{A})=Tr(bcirc(\mathscr{A}))=\sum_i\lambda_i(\mathscr{A})
=p\sum_{i=1}^n(\mathbf{A}^{(1)})_{ii}$;
\item[(c)] $Tr(\mathscr{A}\ast\mathscr{B})=Tr(\mathscr{B}\ast\mathscr{A})$;
\item[(d)] $spec(\mathscr{C}^{-1}\mathscr{A}\mathscr{C})=spec(\mathscr{A})$ and $Tr(\mathscr{C}^{-1}\mathscr{A}\mathscr{C})=Tr(\mathscr{A})$.
\end{itemize}
\end{proposition}

%
%
%
%
%
%

\begin{remark}\label{remark-2}
(i) From Lemma \ref{lemma-1}$(a)$ and Proposition \ref{pro-2}, it is easy to see that for any $\mathscr{A}\in S\mathbb{R}^{n\times n\times p}$ and $\mathscr{B}\in S\mathbb{R}^{n\times n\times p}$,
$$p\langle\mathscr{A},\mathscr{B}\rangle=\langle
bcirc(\mathscr{A}),bcirc(\mathscr{B})\rangle=Tr(bcirc(\mathscr{A})bcirc(\mathscr{B}))=Tr(\mathscr{A}\ast\mathscr{B}),$$
and $\mathscr{A}$ is T-P(S)D if and only if $Tr(\mathscr{V}^\top\ast\mathscr{A}\ast\mathscr{V})> 0$ for all nonzero $\mathscr{V}\in\mathbb{R}^{n\times 1 \times p}$ $(Tr(\mathscr{V}^\top\ast\mathscr{A}\ast\mathscr{V})\geq 0$ for all $\mathscr{V}\in\mathbb{R}^{n\times 1\times p})$.

(ii) Let $\mathscr{A}\in S\mathbb{R}^{n\times n\times p}_{+}$ and $\mathscr{B}\in S\mathbb{R}^{n\times n\times p}_+$. Then, it follows from \cite[Lemmas 1.2.3, 1.2.4]{H-00} and Proposition \ref{pro-2} that
$\langle\mathscr{A},\mathscr{B}\rangle\geq0;\;\; \langle\mathscr{A},\mathscr{B}\rangle=0 \;\mbox{\rm iff}\; \mathscr{A}\ast\mathscr{B}=\mathscr{O};$ and
$$\begin{array}{l}
p\langle\mathscr{A},\mathscr{B}\rangle\geq\lambda_{\min}(\mathscr{A})\lambda_{\max}(\mathscr{B})\leq
\lambda_{\min}(\mathscr{A})Tr(\mathscr{B});\\
p\langle\mathscr{A},\mathscr{B}\rangle\leq
\lambda_{\max}(\mathscr{A})Tr(\mathscr{B})\leq
n\lambda_{\max}(\mathscr{A})\lambda_{\max}(\mathscr{B}).
\end{array}$$
\end{remark}

It is known that the eigenvalue decomposition plays an important role in the study of symmetric matrices. In the following, we establish a similar decomposition for symmetric third-order tensors, especially for the T-P(S)D tensor.
\begin{theorem}(T-eigenvalue decomposition)\label{T-ED}
Every $\mathscr{A}\in S\mathbb{R}^{n\times n\times p}$ can be factored as
$$\mathscr{A}=\mathscr{U}^{\top}\ast\mathscr{S}\ast\mathscr{U},$$
where $\mathscr{U}\in \mathbb{R}^{n\times n\times p}$ is an orthogonal tensor and $\mathscr{S}\in \mathbb{R}^{n\times n\times p}$ is an F-diagonal tensor(That is, each frontal slice of $\mathscr{S}$ is a diagonal matrix) with all
of the diagonal entries of $(\mathbf{F}_{p}\otimes I_{n\times n})bcirc(\mathscr{S})(\mathbf{F}_{p}^H\otimes I_{n\times n})$ being the T-eigenvalues of $\mathscr{A}$. In particular, if $\mathscr{A}\in S\mathbb{R}^{n\times n\times p}_+$ ($\mathscr{A}\in S\mathbb{R}^{n\times n\times p}_{++}$), then all of the diagonal entries of $(\mathbf{F}_{p}\otimes I_{n\times n})bcirc(\mathscr{S})(\mathbf{F}_{p}^H\otimes I_{n\times n})$ are nonnegative (positive).
\end{theorem}

\subsection{The T-roots of a symmetric T-PSD tensor}\label{Sect. 4.4}
The following result about the roots of a symmetric T-PSD tensor also holds by Theorem \ref{pro-1(a)}, Lemma \ref{lemma-1} and \cite[Theorem 7.2.6]{HJ-13}.

\begin{theorem}(The T-roots of a symmetric T-PSD tensor)\label{roots}
Let $\mathscr{A}\in S\mathbb{R}^{n\times n\times p}_+$ and $k\geq 1$. Then there exists a unique $\mathscr{B}\in S\mathbb{R}^{n\times n\times p}_+$ with $\mathscr{B}^k=\mathscr{A}$.
\end{theorem}


\begin{corollary}\label{square-root}
Let $\mathscr{A}\in S\mathbb{R}^{n\times n\times p}$ be T-PSD. Then there exists a unique positive semidefinite tensor $\mathscr{B}\in S\mathbb{R}^{n\times n\times p}$ with $\mathscr{B}^2=\mathscr{A}$. We write such $\mathscr{B}$ as $\mathscr{A}^{\frac{1}{2}}$.
\end{corollary}

Furthermore, the following conclusion is true.
\begin{theorem}\label{equivalent-PSD}
For any $\mathscr{A}\in S\mathbb{R}^{n\times n\times p}$ with $bcirc(\mathscr{A})$ being block diagonalized as (\ref{PSD-A}), (a) $\mathscr{A}\in S\mathbb{R}^{n\times n\times p}_+ (\mathscr{A}\in S\mathbb{R}^{n\times n\times p}_{++})$ if and only if (b) $\mathscr{A}=\mathscr{P}^\top\ast\mathscr{P}$ for some tensor $\mathscr{P}\in\mathbb{R}^{m\times n \times p}$ $(\mathscr{A}=\mathscr{P}^\top\ast\mathscr{P}$ for some nonsingular tensor $\mathscr{P}\in\mathbb{R}^{m\times n \times p}$).
\end{theorem}


\subsection{The cone of T-PSD tensors}\label{Sect. 4.5}
In this subsection, we investigate the set of T-PSD tensors. Recall that a subset $C$ of a vector space $V$ is called a cone (or sometimes called a linear cone) if for each $\mathbf{x}$ in $C$ and any nonnegative scalar $\alpha$, the product $\alpha\mathbf{x}$ is in $C$; $C$ is called a convex cone if for any nonnegative scalars $\alpha,\beta$ and any $\mathbf{x}$, $\mathbf{y}$ in $C$, it follows that $\alpha\mathbf{x}+\beta\mathbf{y}$ belongs to $C$; and if additionally $C$ is a closed set, then we call $C$ a closed, convex cone.

\begin{proposition}
$S\mathbb{R}^{n\times n\times p}$ is isomorphic to $\mathbb{R}^{\frac{pn^2+n}{2}}$ if $p$ is odd; and $S\mathbb{R}^{n\times n\times p}$ is isomorphic to $\mathbb{R}^{\frac{pn^2}{2}+n}$ if $p$ is even.
\end{proposition}

\begin{proof}
This proposition can be easily proved; and we omit the proof here.
\end{proof}

\begin{proposition}\label{prop-cone}
$S\mathbb{R}^{n\times n\times p}_+$ is a nonempty, closed, convex, pointed cone.
\end{proposition}

\begin{proof} From Theorem \ref{pro-1(a)} and the fact that $S\mathbb{R}^{n\times n}$ is nonempty and closed, it follows that $S\mathbb{R}^{n\times n\times p}$ is nonempty and closed. For any $\mathscr{A}, \mathscr{B}\in S\mathbb{R}^{n\times n\times p}$ and any two nonnegative scalars $\alpha$ and $\beta$, we have that $bcirc(\alpha\mathscr{A}+\beta\mathscr{B})=\alpha bcirc(\mathscr{A})+\beta bcirc(\mathscr{B})$. Suppose that $\mathscr{A}$, $\mathscr{B}\in S\mathbb{R}^{n\times n\times p}_+$, then $bcirc(\mathscr{A})$ and $bcirc(\mathscr{B})$ belong to $S\mathbb{R}^{np\times np}_+$ from Theorem \ref{pro-1(a)}. Therefore, from the fact that $S\mathbb{R}^{np\times np}_+$ is a convex cone, it follows that $bcirc(\alpha\mathscr{A}+\beta\mathscr{B})=\alpha bcirc(\mathscr{A})+\beta bcirc(\mathscr{B})\in S\mathbb{R}^{np\times np}_+$ , which together with Theorem \ref{pro-1(a)} implies that $\alpha\mathscr{A}+\beta\mathscr{B}\in S\mathbb{R}^{n\times n\times p}_+$. That is to say, $S\mathbb{R}^{n\times n\times p}_+$ is a convex cone. Suppose that $\mathscr{A}\in S\mathbb{R}^{n\times n\times p}_+$ and $\mathscr{-A}\in S\mathbb{R}^{n\times n\times p}_+$, then $bcirc(\mathscr{A})\in S\mathbb{R}^{np\times np}_+$ and $bcirc(\mathscr{-A})=-bcirc(\mathscr{A})\in S\mathbb{R}^{np\times np}_+$. From the fact that $S\mathbb{R}^{np\times np}_+$ is pointed, it follows that $bcirc(\mathscr{A})=\mathbf{O}$. Thus $\mathscr{A}=\mathscr{O}$, which implies that $S\mathbb{R}^{n\times n\times p}_+$ is pointed.  \end{proof}

\begin{remark}
By the proof of Theorem \ref{prop-cone}, it is easy to obtain that $S\mathbb{R}^{n\times n\times p}_{++}$ is a nonempty, open, convex cone. It is also easy to show that $S\mathbb{R}^{n\times n\times p}_{++}$ is the interior of $S\mathbb{R}^{n\times n\times p}_{+}$.
By the theory of conic optimization, it follows that $S\mathbb{R}^{n\times n\times p}_{+}$ ($S\mathbb{R}^{n\times n\times p}_{++}$) can induce a partial order on $S\mathbb{R}^{n\times n\times p}$, denoted by $\succeq_\mathscr{T}(\succ_\mathscr{T})$. That is, for any $\mathscr{A},\mathscr{B}\in S\mathbb{R}^{n\times n\times p}$, $\mathscr{A}\succeq_\mathscr{T}(\succ_\mathscr{T})\; \mathscr{B}$ if and only if $\mathscr{A}-\mathscr{B}\in S\mathbb{R}^{n\times n\times p}_+(S\mathbb{R}^{n\times n\times p}_{++})$.
\end{remark}

In the following, we will use $\mathscr{A}\succeq_\mathscr{T}(\succ_\mathscr{T})\mathscr{O}$ if $\mathscr{A}\in S\mathbb{R}^{n\times n\times p}_+(S\mathbb{R}^{n\times n\times p}_{++})$. Especially, we replace $\mathscr{A}\succeq_\mathscr{T}(\succ_\mathscr{T})\mathscr{O}$ with $\mathbf{A}\succeq(\succ)\mathbf{O}$ if $\mathbf{A}\in S\mathbb{R}^{n\times n}_+(S\mathbb{R}^{n\times n}_{++})$, as any $\mathscr{A}\in S\mathbb{R}^{n\times n\times p}_+(S\mathbb{R}^{n\times n\times p}_{++})$ reduces to the $\mathbf{A}\in S\mathbb{R}^{n\times n}_+(S\mathbb{R}^{n\times n}_{++})$ when $p=1$.

Recall that $S\mathbb{R}^{n\times n}_+$ is a self-dual cone, which plays an important role in the widely studied semidefinite programming. In the following, we generalize this fact to T-semidefinite cone $S\mathbb{R}^{n\times n\times p}_+$.
For a cone $C$, the polar cone (or dual cone) \cite{H-00} is the set $C^*:=\{\mathbf{y}:\langle\mathbf{x},\mathbf{y}\rangle\geq 0, \;\; \mbox{\rm for\ any\ } \mathbf{x}\in C\}$.

\begin{theorem}(Self-duality)\label{self-dual}
$S\mathbb{R}^{n\times n\times p}_+=(S\mathbb{R}^{n\times n\times p}_+)^*$.
\end{theorem}

\begin{proof} (i) $S\mathbb{R}^{n\times n\times p}_+\subseteq(S\mathbb{R}^{n\times n\times p}_+)^*$: To this end, we only need to show that $\mathscr{A}\bullet\mathscr{B}\geq 0$ for any $\mathscr{A}$, $\mathscr{B}\succeq_\mathscr{T}\mathscr{O}$. Since $\mathscr{A}$, $\mathscr{B}\succeq_\mathscr{T}\mathscr{O}$, it follows from Theorem \ref{square-root} that there exists $\mathscr{A}^{\frac{1}{2}}$ and $\mathscr{B}^{\frac{1}{2}}$ such that $\mathscr{A}=\mathscr{A}^{\frac{1}{2}}\ast\mathscr{A}^{\frac{1}{2}}$ and $\mathscr{B}=\mathscr{B}^{\frac{1}{2}}\ast\mathscr{B}^{\frac{1}{2}}$. Thus, we can obtain that $$\begin{array}{lcl}
\mathscr{A}\bullet\mathscr{B}&=&\frac{1}{p}Tr(\mathscr{A}\ast\mathscr{B})=
\frac{1}{p}Tr(\mathscr{A}^{\frac{1}{2}}\ast\mathscr{A}^{\frac{1}{2}}\ast\mathscr{B}^{\frac{1}{2}}\ast\mathscr{B}^{\frac{1}{2}})\\
&=&\frac{1}{p}Tr(\mathscr{B}^{\frac{1}{2}}\ast\mathscr{A}^{\frac{1}{2}}\ast\mathscr{A}^{\frac{1}{2}}\ast\mathscr{B}^{\frac{1}{2}})
=(\mathscr{A}^{\frac{1}{2}}\ast\mathscr{B}^{\frac{1}{2}})\bullet(\mathscr{A}^{\frac{1}{2}}\ast\mathscr{B}^{\frac{1}{2}})\geq 0.
\end{array}$$

(ii) $(S\mathbb{R}^{n\times n\times p}_+)^*\subseteq S\mathbb{R}^{n\times n\times p}_+$: We only need to show $\mathscr{A}\succeq_\mathscr{T}\mathscr{O}$ if $\mathscr{A}\in (S\mathbb{R}^{n\times n\times p}_+)^*$. Suppose $\mathscr{A}\in (S\mathbb{R}^{n\times n\times p}_+)^*$, then $\mathscr{A}\bullet\mathscr{B}\geq 0$ for any $\mathscr{B}\succeq_\mathscr{T}\mathscr{O}$. Taken $\mathscr{B}=\mathscr{D}\ast\mathscr{D}^\top$ where $\mathscr{D}\in \mathbb{R}^{n\times 1\times p}$ is an arbitrary given tensor, then we have that $\mathscr{B}\succeq_\mathscr{T}\mathscr{O}$ and $\mathscr{A}\bullet(\mathscr{D}\ast\mathscr{D}^\top)\geq 0$, i.e., $Tr(\mathscr{A}\ast(\mathscr{D}\ast\mathscr{D}^\top))\geq 0$. By Proposition \ref{pro-2}$(b)$, we get  $$\begin{array}{lcl}
Tr(\mathscr{D}^\top\ast\mathscr{A}\ast\mathscr{D})&=&Tr(bcirc(\mathscr{D}^\top)bcirc(\mathscr{A})bcirc(\mathscr{D}))\\
&=&Tr(bcirc(\mathscr{A})(bcirc(\mathscr{D})bcirc(\mathscr{D}^\top)))=Tr(\mathscr{A}\ast(\mathscr{D}\ast\mathscr{D}^\top))\geq 0.\end{array}$$
Hence, by Remark \ref{remark-2}(i), we can obtain that $\mathscr{A}\succeq_\mathscr{T}\mathscr{O}$.
\end{proof}

\subsection{The T-schur complement of a symmetric T-PSD tensor}\label{Sect. 4.6}
In this subsection, we give a characterization of the T-positive semidefiniteness of a third-order tensor by the T-positive semidefiniteness of the T-schur complement \cite{MQW-19}.

\begin{lemma}(Tensor block multiplication via T-product)\label{Block-Multiplication} {\rm\cite{MQW-19-1}}
Suppose $\mathscr{A}_1\in \mathbb{C}^{n_1\times m_1\times p}$, $\mathscr{B}_1\in \mathbb{C}^{n_1\times m_2\times p}$, $\mathscr{C}_1\in \mathbb{C}^{n_2\times m_1\times p}$, $\mathscr{D}_1\in \mathbb{C}^{n_2\times m_2\times p}$, $\mathscr{A}_2\in \mathbb{C}^{m_1\times r_1\times p}$, $\mathscr{B}_2\in \mathbb{C}^{m_1\times r_2\times p}$, $\mathscr{C}_2\in \mathbb{C}^{m_2\times r_1\times p}$ and $\mathscr{D}_2\in
\mathbb{C}^{m_2\times r_2\times p}$ are complex tensors, then
$$\left[
    \begin{array}{cc}
      \mathscr{A}_1 & \mathscr{B}_1 \\
      \mathscr{C}_1 & \mathscr{D}_1 \\
    \end{array}
  \right]
\ast \left[
    \begin{array}{cc}
      \mathscr{A}_2 & \mathscr{B}_2 \\
      \mathscr{C}_2 & \mathscr{D}_2 \\
    \end{array}
  \right]=\left[
    \begin{array}{cc}\mathscr{A}_1\ast\mathscr{A}_2+\mathscr{B}_1\ast \mathscr{C}_2 &\mathscr{A}_1\ast \mathscr{B}_2+ \mathscr{B}_1\ast \mathscr{D}_2\\
\mathscr{C}_1\ast\mathscr{A}_2 +\mathscr{D}_1\ast\mathscr{C}_2 & \mathscr{C}_1\ast\mathscr{B}_2 +\mathscr{D}_1\ast\mathscr{D}_2 \\
    \end{array}
  \right].$$
\end{lemma}

\begin{lemma}\label{T-block-PSD}
Suppose that $\mathscr{A}_i\in S\mathbb{R}^{n_i\times n_i\times p}$ for any $i\in [m]$. Then the block diagonal tensor
$\mathscr{A}=Diag(\mathscr{A}_i: i\in [p])$ is symmetric T-P(S)D if and only if all $\mathscr{A}_i$ are so.
\end{lemma}

\begin{proof} For any nonzero $\mathscr{V}\in\mathbb{R}^{n\times 1 \times p}$ with $n=\sum_{i=1}^{m}n_i$, we divided it into a block tensor, i.e., $\mathscr{V}=vec(\mathscr{V}_i: i\in [m])$
where $\mathscr{V}_{i}\in \mathbb{R}^{n_i\times 1 \times p}$ for any $i\in [m]$.
Then, by Lemma \ref{Block-Multiplication}, we can obtain that
$$\begin{array}{lcl}
\mathscr{V}^\top\ast \mathscr{A}\ast \mathscr{V}
&=&vec(\mathscr{V}_i: i\in [m])^\top\ast Diag(\mathscr{A}_i: i\in [p])\ast vec(\mathscr{V}_i: i\in [m])\\
&=&\sum_{i=1}^m\mathscr{V}_{i}^\top\ast\mathscr{A}_{i}\ast\mathscr{V}_{i}.
\end{array}$$
Hence, it is not difficult to get that $\mathscr{A}$ is symmetric T-P(S)D iff all $\mathscr{A}_i$ are so.
\end{proof}

Besides, from \cite[Proposition 1.1.7]{H-00}, Lemma \ref{lemma-1} and Theorem \ref{pro-1(a)}, the following result holds.
\begin{lemma}\label{coro-1}
Suppose that $\mathscr{B}\in \mathbb{R}^{n\times n\times p}$ be nonsingular. Then $\mathscr{A}\in S\mathbb{R}^{n\times n\times p}_{+}$ ($S\mathbb{R}^{n\times n\times p}_{++}$) if and only if $(\mathscr{B}^\top\ast\mathscr{A}\ast\mathscr{B})\in S\mathbb{R}^{n\times n\times p}_{+}$ ($S\mathbb{R}^{n\times n\times p}_{++}$).
\end{lemma}


Now, we can establish a theorem about the T-Schur complement.

\begin{theorem}\label{T-schur}(T-Schur complement)
Suppose that $\mathscr{A}\in S\mathbb{R}^{m\times m\times p}_{++}$, $\mathscr{C}\in S\mathbb{R}^{n\times n\times p}$, and $\mathscr{B}\in \mathbb{R}^{m\times n\times p}$. Then
$$\begin{array}{ccc}
\begin{bmatrix}
  \mathscr{A} & \mathscr{B} \\
  \mathscr{B}^\top & \mathscr{C}
\end{bmatrix}\succ_\mathscr{T}(\succeq_\mathscr{T})\; \mathscr{O}&\Longleftrightarrow& \mathscr{C}-\mathscr{B}^\top\ast\mathscr{A}^{-1}\ast\mathscr{B}\succ_\mathscr{T} (\succeq_\mathscr{T})\;\mathscr{O}
\end{array}.$$
\end{theorem}

\begin{proof} It follows from $\mathscr{A}\in S\mathbb{R}^{m\times m\times p}_{++}$ that $\mathscr{A}$ is nonsingular. Denote the block tensor
$$\mathscr{D}:=\begin{bmatrix}
  \mathscr{I}_{mmp} & -\mathscr{A}^{-1}\ast\mathscr{B} \\
  \mathscr{O} & \mathscr{I}_{nnp}
\end{bmatrix},$$
then we have
$$\begin{array}{ccc}
&\mathscr{D}^\top\ast \begin{bmatrix}
  \mathscr{A} & \mathscr{B} \\
  \mathscr{B}^\top & \mathscr{C}
\end{bmatrix}\ast \mathscr{D}&
=\begin{bmatrix}
  \mathscr{A} & \mathscr{O} \\
  \mathscr{O} & \mathscr{C}-\mathscr{B}^\top\ast\mathscr{A}^{-1}\ast\mathscr{B}
\end{bmatrix}.
\end{array}$$
Therefore, by Lemma \ref{T-block-PSD} and Lemma \ref{coro-1}, the theorem is proved. \end{proof}

\section{Semidefinite programming over the third-order symmetric tensor space}
\label{Sect. 4}

In this section, we first introduce the TSDP and give its duality theory; and then, we show the transformation of TSDPs into SDPs in the complex domain. After that, we consider several problems and reformulate (or relax) them as TSDPs. Finally, we present some preliminary numerical results for solving the unconstrained polynomial optimization problem via the TSDP relaxation.

\subsection{ TSDP problems in primal-dual forms}\label{Sect. 5.1}
In this subsection, we replace the matrix variables in the classic SDP by the tensor variables to yield the TSDP. We consider the TSDP in primal form:

$$
\mbox{\rm (PTSDP)}\quad \min\limits_{\mathscr{X}}\;  \langle\mathscr{C}, \mathscr{X}\rangle\quad
\mbox{\rm s.t.}\quad \mathscr{A}\mathscr{X}=[\langle\mathscr{A}_i,\mathscr{X}\rangle]_{i\in [m]}=\mathbf{b},
\;\; \mathscr{X}\succeq_\mathscr{T} \mathscr{O},
$$
where all $\mathscr{A}_i\in S\mathbb{R}^{n\times n\times p}$, $\mathbf{b}\in \mathbb{R}^{m}$, $\mathscr{C}\in S\mathbb{R}^{n\times n\times p}$ are given and $\mathscr{X}\in S\mathbb{R}^{n\times n\times p}$ is the variable. $\mathscr{A}$ is a linear operator from $S\mathbb{R}^{n\times n\times p}$ into $\mathbb{R}^{m}$.

Just as the derivation of the dual problem of the SDP, in order to obtain the dual problem of (PTSDP), we try to find the adjoint operator of $\mathscr{A}$ at first, which is a linear operator from $\mathbb{R}^{m}$ into $S\mathbb{R}^{n\times n\times p}$ satisfying $\langle\mathscr{A}\mathscr{X},\mathbf{y}\rangle=\langle\mathscr{X},\mathscr{A}^*\mathbf{y}\rangle$ for any $\mathscr{X}$ in $S\mathbb{R}^{n\times n\times p}$ and $\mathbf{y}$ in $\mathbb{R}^{m}$. Since
$$\langle\mathscr{A}\mathscr{X},\mathbf{y}\rangle=\frac{1}{p}\sum_{i=1}^m {y}_iTr(\mathscr{A}_i\ast\mathscr{X})=\frac{1}{p}Tr(\mathscr{X}\ast\sum_{i=1}^m{y}_i\mathscr{A}_i)=\langle\mathscr{X},\mathscr{A}^*\mathbf{y}\rangle,$$
we have $\mathscr{A}^*\mathbf{y}=\sum_{i=1}^m{y}_i\mathscr{A}_i$. Now we can construct the dual of (PTSDP) by the Lagrange approach. By adding a Lagrange multiplier $\mathbf{y\in \mathbb{R}}^{m}$, (PTSDP) can be turned into $\inf_{\mathscr{X}\succeq_\mathscr{T} \mathscr{O}}\sup_{\mathbf{y\in \mathbb{R}}^{m}}\langle\mathscr{C}, \mathscr{X}\rangle+\langle\mathbf{b}-\mathscr{A}\mathscr{X},\mathbf{y}\rangle$, then the dual of (PTSDP) yields through interchanging $\inf$ and $\sup$.
Note that
$$\sup_{\mathbf{y}\in \mathbb{R}^{m}}\inf_{\mathscr{X}\succeq_\mathscr{T} \mathscr{O}}\langle\mathbf{b}, \mathbf{y}\rangle+\langle\mathscr{C}-\mathscr{A}^*\mathbf{y},\mathscr{X}\rangle=
\left\{\begin{array}{ll}
 \langle\mathbf{b},\mathbf{y}\rangle, & \mbox{\rm if}\;\;\mathscr{C}-\mathscr{A}^*\mathbf{y}\in (S\mathbb{R}^{n\times n\times p}_+)^*,\\
 -\infty, & \mbox{\rm otherwise}.
\end{array}\right.
$$
This, together with Theorem \ref{self-dual}, implies that we can write the dual problem of (PTSDP) by introducing a slack tensor $\mathscr{S}$ as:
$$
\mbox{\rm(DTSDP)}\quad \max\limits_{\mathbf{y},\mathscr{S}}\; \langle\mathbf{b},\mathbf{y}\rangle\quad
 \mbox{\rm s.t.}\quad \mathscr{A}^*\mathbf{y}+\mathscr{S}=\mathscr{C},
\;\; \mathscr{S}\succeq_\mathscr{T} \mathscr{O},
$$
where $\mathbf{y}\in \mathbb{R}^{m}$ and $\mathscr{S}\in S\mathbb{R}^{n\times n\times p}$ are the variables. When $p=1$, the TSDP problems (PTSDP) and (DTSDP) are corresponding to the classic SDP problems in primal-dual forms. Denote
\begin{eqnarray}\label{F-D}\begin{array}{l}
F(P):=\{\mathscr{X}\in S\mathbb{R}^{n\times n \times p}: \mathscr{A}\mathscr{X}=\mathbf{b},\;\mathscr{X}\succeq_\mathscr{T} \mathscr{O}\},\\
F(D):= \{(\mathbf{y},\mathscr{S})\in \mathbb{R}^m \times S\mathbb{R}^{n\times n\times p}: \mathscr{A}^*\mathbf{y}+\mathscr{S}=\mathscr{C}, \; \mathscr{S}\succeq_\mathscr{T} \mathscr{O}\},\\
p^*:=\inf\{\langle\mathscr{C},{\mathscr{X}}\rangle: \mathscr{X}\in F(P)\}\;\;\mbox{\rm and}\;\;
d^*:=\sup\{\langle\mathbf{b},\mathbf{y}\rangle: (\mathbf{y},\mathscr{S})\in F(D)\}.
\end{array}\end{eqnarray}

From properties of the T-semidefinite cone obtained in Section \ref{Sect. 4.5} and the theory of conic optimization problems \cite{BN-01}, it is not difficult to obtain the following
results, and the proofs are omitted here.

\begin{theorem}
Let $F(P)$, $F(D)$, $p^*$ and $d^*$ be defined as (\ref{F-D}). Suppose that $\mathscr{X}\in F(P)$ and $(\mathbf{y},\mathscr{S})\in F(D)$. Then
\begin{itemize}
\item{(weak duality)\;\;$\langle b,y\rangle\leq\langle\mathscr{C},\mathscr{X}\rangle$.}
\item{(strong duality)\;\;Suppose that \mbox{\rm (PTSDP)} is bounded below and strictly feasible (respectively, \mbox{\rm (DTSDP)} is
bounded above and strictly feasible), then $p^*=d^*$ and \mbox{\rm(DTSDP)} (respectively, \mbox{\rm (PTSDP)}) is solvable.}
\item{(complementarity slackness condition)\;\;If $p^*=d^*$, then $\mathscr{X}$ is optimal for \mbox{\rm(PTSDP)} and $(\mathbf{y},\mathscr{S})$  is optimal for \mbox{\rm(DTSDP)} if and only if the complementarity slackness condition holds, that is, $\langle\mathscr{X},\mathscr{S}\rangle=0$.}
\item{(optimality condition)\;\;If $\langle\mathscr{C},\mathscr{X}\rangle=\langle b,y\rangle$, then $\mathscr{X}$ is optimal for \mbox{\rm (PTSDP)}, and $(\mathbf{y},\mathscr{S})$ is optimal for \mbox{\rm (DTSDP)}.}
\end{itemize}

\end{theorem}

\subsection{The transformation of TSDPs into SDPs in the complex domain}\label{Sect. 5.2}
In this subsection, we present a method to solve the TSDP problem by transforming it into an SDP in the complex domain (CSDP for short).

For any $\mathscr{A}_i\in S\mathbb{R}^{n\times n \times p}$ ($i\in [m]$) in (PTSDP), $bcirc(\mathscr{A}_i)$ can be block diagonalized as
$bcirc({\mathscr{A}_i})=(\mathbf{F}^{H}_{p}\otimes\mathbf{I}_{n\times n})\mathbf{A}^i(\mathbf{F}_{p}\otimes\mathbf{I}_{n\times n})$
with $\mathbf{A}^i=Diag(\mathbf{A}^i_j: j\in [p])$
where all $\mathbf{A}^i_j\in H\mathbb{C}^{n\times n}$ ($\mathbf{A}^i_j\in H\mathbb{C}^{n\times n}_+$ if particularly $\mathscr{A}\in S\mathbb{R}^{n\times n \times p}_+$). Note that
$$\begin{array}{lcl}
&&\langle\mathscr{C},{\mathscr{X}^*}\rangle=\min_{\mathscr{X}}\langle\mathscr{C},{\mathscr{X}}\rangle\\
&\Leftrightarrow& \langle bcirc(\mathscr{C}),bcirc({\mathscr{X}^*})\rangle=\min_{bcirc(\mathscr{X})}\langle bcirc(\mathscr{C}),bcirc({\mathscr{X}})\rangle\\
&\Leftrightarrow& \langle\mathbf{C},{\mathbf{X}^*}\rangle=\min_{\mathbf{X}}\langle\mathbf{C},\mathbf{X}\rangle,
\end{array}$$
where $bcirc(\mathscr{C})=(\mathbf{F}_{p}^{H}\otimes\mathbf{I}_{n\times n})\mathbf{C}(\mathbf{F}_{p}\otimes\mathbf{I}_{n\times n})$, $bcirc(\mathscr{X})=(\mathbf{F}_{p}^{H}\otimes\mathbf{I}_{n\times n})\mathbf{X}(\mathbf{F}_{p}\otimes\mathbf{I}_{n\times n})$ and $bcirc({\mathscr{X}^*})=(\mathbf{F}_{p}^{H}\otimes\mathbf{I}_{n\times n}){\mathbf{X}^*}(\mathbf{F}_{p}\otimes\mathbf{I}_{n\times n})$ with
\begin{eqnarray}\label{S}
\mathbf{C}=Diag(\mathbf{C}_i,i\in[p]),
\mathbf{X}=Diag(\mathbf{X}_i,i\in[p]),
\mathbf{X}^*=Diag(\mathbf{X}^*_i,i\in[p])
\end{eqnarray}
with all $\mathbf{C}_i$, $\mathbf{X}_i$ and ${\mathbf{X}^*}_i$ in $H\mathbb{C}^{n\times n}$.
In addition, $\mathscr{X}\succeq_\mathscr{T}\mathscr{O}\Leftrightarrow \mathbf{X}\succeq\mathbf{O}$ and $$\mathscr{A}\mathscr{X}=[\langle\mathscr{A}_i, \mathscr{X}\rangle]_{i\in [m]}=
[\frac{1}{p}Tr(bcirc(\mathscr{A}_i),bcirc(\mathscr{X}))]_{i\in [m]}=
[\frac{1}{p}\langle\mathbf{A}^i, \mathbf{X}\rangle]_{i\in [m]}.
$$
 Therefore, let $S$ denote the space of block diagonal Hermitian matrices as the form in (\ref{S}),
then (PTSDP) and (DTSDP) are equivalent to the following SDP problems (PCSDP) and (DCSDP) respectively:
$$\begin{array}{lcl}
\mbox{\rm(PCSDP)}\quad & \min\limits_{\mathbf{X}\in S}\; \frac{1}{p}\langle\mathbf{C}, \mathbf{X}\rangle\quad
 \mbox{\rm s.t.}\quad \mathbf{A}\mathbf{X}=p\mathbf{b},\; \mathbf{X}\succeq \mathbf{O},\vspace{2mm}\\
\mbox{\rm(DCSDP)}\quad & \max\limits_{(\mathbf{y},\mathbf{S})\in \mathbb{R}^m\times S}\; \langle\mathbf{b},\mathbf{y}\rangle\quad
\mbox{\rm s.t.}\quad \mathbf{A}^*\mathbf{y}+\mathbf{S}=\mathbf{C},\;
 \mathbf{S}\succeq \mathbf{O},
\end{array}$$
where $\mathbf{X}$, $\mathbf{C}$ and $\mathbf{S}$ are given as (\ref{S}), and $\mathbf{A}$ is a linear operator from $H\mathbb{C}^{np\times np}$ into $\mathbb{R}^{m}$ denoted as $\mathbf{A}\mathbf{X}=[\langle\mathbf{A}^i, \mathbf{X}\rangle]_{i\in[m]}$ with $\mathbf{A}^*$ being its adjoint operator.

It should be noted that both (DCSDP) and (PCSDP) are SDPs in the complex domain and (DCSDP) is exactly the dual problem of (PCSDP). Noting that these diagonal blocks of the complex matrices $\mathbf{X}$, $\mathbf{C}$ and $\mathbf{A}^i$ for $i\in [m]$, obtained by block diagonalizing the real tensors $\mathscr{X}$, $\mathscr{C}$ and $\mathscr{A}^i$, satisfy the relationships described in Lemma \ref{coro-*}. So, (PCSDP) and (DCSDP) can be converted to SDPs of smaller size. For the cleanness  of the paper, we only take the transformation of (PCSDP) for example, which can be divided into the following two cases.

{\bf Case 1: $p$ is even.} Let $\mathbf{X}$, $\mathbf{C}$ and $\mathbf{A}^i$ for $i\in [m]$ be the complex block diagonal matrices in (PCSDP), which are obtained by block diagonalizing the real tensors $\mathscr{X}$, $\mathscr{C}$ and $\mathscr{A}^i$ for $i\in [m]$ in (PTSDP), respectively. From Lemma \ref{coro-*}, it follows that for any $j\in [p]\setminus\{1,\frac{p+2}{2}\}$ and $i\in [m]$,
$$\left\{\begin{array}{lll}
\mathbf{X}_1\in \mathbb{R}^{n\times n},& \mathbf{X}_{\frac{p+2}{2}}\in \mathbb{R}^{n\times n},& \mathbf{X}_j\in \mathbb{C}^{n\times n},\;\; \mathbf{X}_j=\overline{\mathbf{X}_{p+2-j}};\\
\mathbf{C}_1\in \mathbb{R}^{n\times n},& \mathbf{C}_{\frac{p+2}{2}}\in \mathbb{R}^{n\times n},& \mathbf{C}_j\in \mathbb{C}^{n\times n},\;\; \mathbf{C}_j=\overline{\mathbf{C}_{p+2-j}};\\
\mathbf{A}^i_1\in \mathbb{R}^{n\times n},& \mathbf{A}_{\frac{p+2}{2}}\in \mathbb{R}^{n\times n},& \mathbf{A}^i_j\in \mathbb{C}^{n\times n},\;\; \mathbf{A}^i_j=\overline{\mathbf{A}_{p+2-j}}.\\
\end{array}\right.$$
Thus, for any $i\in [m]$,
$$\begin{array}{rcl}
\mathbf{A}^i\bullet\mathbf{X}&=&\mathbf{A}^i_1\bullet\mathbf{X}_1+\mathbf{A}^i_2\bullet\mathbf{X}_2+
\cdots+\mathbf{A}^i_p\bullet\mathbf{X}_p\\
 &=&\mathbf{A}^i_1\bullet\mathbf{X}_1+\mathbf{A}^i_{\frac{p+2}{2}}\bullet\mathbf{X}_{\frac{p+2}{2}}+
 \sum_{j=1}^\frac{p}{2}(\mathbf{A}^i_j\bullet\mathbf{X}_j+\overline{\mathbf{A}^i_j}\bullet\overline{\mathbf{X}_j})\\
 &=&\mathbf{A}^i_1\bullet\mathbf{X}_1+\mathbf{A}^i_{\frac{p+2}{2}}\bullet\mathbf{X}_{\frac{p+2}{2}}+
 2\sum_{j=1}^\frac{p}{2}\mathbf{A}^i_j\bullet\mathbf{X}_j,
\end{array}$$
where the last equality follows from the fact that the inner product between two Hermitian matrices is real. Similarly, we can also obtain that $$\mathbf{C}\bullet\mathbf{X}=\mathbf{C}_1\bullet\mathbf{X}_1 +\mathbf{C}_{\frac{p+2}{2}}\bullet\mathbf{X}_{\frac{p+2}{2}}+
2\sum_{j=1}^\frac{p}{2}\mathbf{C}_j\bullet\mathbf{X}_j.$$ Thus, by letting  $\widetilde{\mathbf{X}}=Diag(\mathbf{X}_1,\mathbf{X}_2,\cdots, \mathbf{X}_\frac{p}{2},\mathbf{X}_\frac{p+2}{2})$, $\widetilde{\mathbf{A}^i}=Diag(\mathbf{A}^i_1,2\mathbf{A}^i_2,\cdots, 2\mathbf{A}^i_\frac{p}{2},\mathbf{A}^i_\frac{p+2}{2})$
for any $i\in [m]$, and $\widetilde{\mathbf{C}}=Diag(\mathbf{C}_1,2\mathbf{C}_2,\cdots, 2\mathbf{C}_\frac{p}{2}$, $\mathbf{C}_\frac{p+2}{2})$,
it follows that (PCSDP) is equivalent to
$$
(\mbox{\rm P}'\mbox{\rm CSDP})\quad \min\limits_{\widetilde{\mathbf{X}}\in S}\; \frac{1}{p}\langle\widetilde{\mathbf{C}}, \widetilde{\mathbf{X}}\rangle\quad
\mbox{\rm s.t.}\quad \widetilde{\mathbf{A}}\widetilde{\mathbf{X}}=p\mathbf{b},\;\;\widetilde{\mathbf{X}}\succeq \mathbf{O},
$$
where $\widetilde{\mathbf{A}}$ is a linear operator with $\widetilde{\mathbf{A}}\widetilde{\mathbf{X}}=[\langle\widetilde{\mathbf{A}^i},\widetilde{\mathbf{X}}\rangle]_{i\in[m]}$.

{\bf Case 2: $p$ is odd.} By the same process as {\bf{Case 1}}, it is not difficult to obtain that (PCSDP) is equivalent to
$$
(\mbox{\rm P}''\mbox{\rm CSDP})\quad \min\limits_{\widetilde{\mathbf{X}}\in S}\; \frac{1}{p}\langle\widetilde{\mathbf{C}}, \widetilde{\mathbf{X}}\rangle\quad
\mbox{\rm s.t.}\quad \widetilde{\mathbf{A}} \widetilde{\mathbf{X}}=p\mathbf{b},\;\;\widetilde{\mathbf{X}}\succeq \mathbf{O},
$$
where $\widetilde{\mathbf{A}^i}=Diag(\mathbf{A}^i_1,2\mathbf{A}^i_2,\cdots, 2\mathbf{A}^i_\frac{p+1}{2})$ for any $i\in [m]$; $\widetilde{\mathbf{C}}=Diag(\mathbf{C}_1,2\mathbf{C}_2,\cdots,2\mathbf{C}_\frac{p+1}{2})$;
$\widetilde{\mathbf{X}}=Diag(\mathbf{X}_1,\mathbf{X}_2,\cdots,\mathbf{X}_\frac{p+1}{2})$;
and $\widetilde{\mathbf{A}}$ is a linear operator with $\widetilde{\mathbf{A}} \ \widetilde{\mathbf{X}}=[\langle\widetilde{\mathbf{A}^i},\widetilde{\mathbf{X}}\rangle]_{i\in[m]}$.

As can be seen from the above discussion, we provide a way to deal with (PTSDP) of size $n\times n\times p$ by transforming it into a CSDP with block diagonal structure of size $n(\frac{p+1}{2})\times n(\frac{p+1}{2})$ as (P$'$CSDP) or $n(\frac{p+2}{2})\times n(\frac{p+2}{2})$ as (P$''$TSDP), which are almost half the size of (PCSDP) with block diagonal structure of size $np\times np$ when $p>2$.

\subsection{Some applications of TSDPs}\label{Sect. 5.3}
In this subsection, we show several applications which can be formulated as TSDP problems.

{\bf Application 1. Minimizing the maximum T-eigenvalue of a third-order symmetric tensor.} The T-eigenvalue was first proposed for third-order symmetric tensors in \cite{MQW-19} by Miao, Qi and Wei. Suppose that $\mathscr{M}(\mathbf{z})\in S\mathbb{R}^{n\times n\times p}$ is a third-order symmetric tensor, which depends linearly on a vector $\mathbf{z}$. Since $\lambda_{max}(\mathscr{M}(\mathbf{z}))\leq\eta$ if and only if $\lambda_{max}(\mathscr{M}(\mathbf{z})-\eta\mathscr{I}_{nnp})\leq 0$, i.e.,  $\lambda_{min}(\eta\mathscr{I}_{nnp}-\mathscr{M}(\mathbf{z}))\geq 0$, which and Remark \ref{remark-3} imply that $\eta\mathscr{I}_{nnp}-\mathscr{M}(\mathbf{z})\succeq_\mathscr{T} \mathscr{O}$. Therefore, the problem of minimizing the maximum T-eigenvalue of $\mathscr{M}(\mathbf{z})$ can be transformed as the following TSDP problem:
\begin{eqnarray*}
\max\limits_{\eta,\mathbf{z}}\; -\eta\quad
\mbox{\rm s.t.}\quad \eta\mathscr{I}_{nnp}-\mathscr{M}(\mathbf{z})\succeq_\mathscr{T} \mathscr{O}.
\end{eqnarray*}

{\bf Application 2. Minimizing the spectral norm of a third-order tensor.} Recall that for any $\mathscr{A}\in \mathbb{R}^{m\times n\times p}$, the tensor spectral norm $\|\mathscr{A}\|_2$ of $\mathscr{A}$ is defined as the largest singular value of $\mathscr{A}$ (see \cite{Z-15,LFCL-19}). It is known that the tensor spectral norm plays an important role in the proof of the optimal conditions for the relative problems in \cite{Z-15,LFCL-19}.

Suppose that $\mathscr{P}(\mathbf{z})\in \mathbb{R}^{m\times n\times p}$ is a third-order real tensor, which depends linearly on a vector $\mathbf{z}$. Noting that $\eta\geq \|\mathscr{P}(\mathbf{z})\|_2$ if and only if $\eta^2\geq\lambda_{max}(\mathscr{P}(\mathbf{z})^\top\ast\mathscr{P}(\mathbf{z}))$; and by Theorem \ref{T-schur}, the latter is equivalent to
$$
\left[\begin{array}{cc}
\eta\mathscr{I}_{mmp} & \mathscr{P}(\mathbf{z})\\
\mathscr{P}(\mathbf{z})^\top & \eta\mathscr{I}_{nnp}
\end{array}\right]\succeq_\mathscr{T} \mathscr{O}.
$$
Therefore, the problem of minimizing $\|\mathscr{P}(\mathbf{z})\|_2$  can be transformed as the following TSDP problem:
\begin{eqnarray*}
\max\limits_{\eta,\mathbf{z}}\; -\eta\quad
\mbox{\rm s.t.}\quad \left[\begin{array}{ll}
\eta\mathscr{I}_{mmp} & \mathscr{P}(\mathbf{z})\\
\mathscr{P}(\mathbf{z})^\top & \eta\mathscr{I}_{nnp}
\end{array}\right]\succeq_\mathscr{T} \mathscr{O}.
\end{eqnarray*}

{{\bf Application 3. Minimizing the nuclear norm of a third-order tensor.}
Recall that the tensor nuclear norm $\|\mathscr{A}\|_*$ of any $\mathscr{A}\in \mathbb{R}^{m\times n\times p}$ is defined as the sum of singular values of the first frontal slice and is shown to be the dual norm of the tensor spectral norm $\|\mathscr{A}\|_2$ in \cite{LFCL-19}. In this part, we consider the following two cases.

\noindent {\bf Case 1. Minimizing the nuclear norm of a third-order tensor without constraint.}
Recall that for a given norm $\|\cdot\|$ in the inner product space consisting of three-order tensors, its dual norm $\|\cdot\|_d$ is defined as
$\|\mathscr{X}\|_d:=\sup\{\langle\mathscr{X},\mathscr{Y}\rangle: \|\mathscr{Y}\|\leq1\}.$
Since the nuclear norm $\|\mathscr{A}\|_*$ of tensor $\mathscr{A}$ is the dual norm of the tensor spectral norm $\|\mathscr{A}\|_2$, according to the relationship that for any $\mathscr{Y}\in \mathbb{R}^{m\times n\times p}$,
$$\begin{array}{rcl}
\|\mathscr{Y}\|_2\leq1&\Longleftrightarrow& \left[\begin{array}{ll}
\mathscr{I}_{mmp} & \mathscr{X}\\
\mathscr{X}^\top &\mathscr{I}_{nnp}
\end{array}\right]\succeq_\mathscr{T} \mathscr{O}
\end{array}$$
as shown in {\bf Application 2} by taking $\eta=1$, we can obtain that the problem of computing the nuclear norm of $\mathscr{A}\in\mathbb{R}^{m\times n\times p}$ is equivalent to the following TSDP model:
\begin{eqnarray}\label{nuclear-norm-TSDP}
\\ \max\limits_\mathscr{X}\; \langle\mathscr{A},\mathscr{X}\rangle=\frac{1}{2}\langle\left[\begin{array}{ll}
\mathscr{O} & \mathscr{A}\\
\mathscr{A}^\top &\mathscr{O}
\end{array}\right]
,\left[\begin{array}{ll}
\mathscr{I}_{mmp} & \mathscr{X}\\
\mathscr{X}^\top &\mathscr{I}_{nnp}
\end{array}\right]\rangle \quad
\mbox{\rm s.t.}\quad \left[\begin{array}{ll}
\mathscr{I}_{mmp} & \mathscr{X}\\
\mathscr{X}^\top &\mathscr{I}_{nnp}
\end{array}\right]\succeq_\mathscr{T} \mathscr{O}.\nonumber
\end{eqnarray}
Noting that for any $\mathscr{W}_1\in \mathbb{R}^{m\times m\times p}$ and $\mathscr{W}_2\in \mathbb{R}^{n\times n\times p}$,
$$\langle\left[\begin{array}{ll}
\mathscr{W}_1 & \mathscr{O}\\
\mathscr{O}^\top &\mathscr{W}_2
\end{array}\right],\left[\begin{array}{ll}
\mathscr{I}_{mmp} & \mathscr{X}\\
\mathscr{X}^\top &\mathscr{I}_{nnp}
\end{array}\right]\rangle=\frac{1}{p}[Tr(\mathscr{W}_1)+Tr(\mathscr{W}_2)],$$
thus, the dual problem of (\ref{nuclear-norm-TSDP}) can be formulated as
\begin{eqnarray}\label{nuclear-norm-TSDP-dual}
\min\limits_{\mathscr{W}_1,\mathscr{W}_2}\; \frac{1}{2p}[Tr(\mathscr{W}_1)+Tr(\mathscr{W}_2)] \quad
\mbox{\rm s.t.}\quad \left[\begin{array}{ll}
\mathscr{W}_1 & \mathscr{A}\\
\mathscr{A}^\top &\mathscr{W}_2
\end{array}\right]\succeq_\mathscr{T} \mathscr{O}.
\end{eqnarray}

It is not difficult to show that there is no duality gap between (\ref{nuclear-norm-TSDP}) and (\ref{nuclear-norm-TSDP-dual}). We omit the proof here.

\noindent {\bf Case 2. Minimizing the nuclear norm of a third-order tensor with an equality constraint.}
In this part, we investigate the problem of minimizing the nuclear norm of $\mathscr{X}\in \mathbb{R}^{m\times n\times p}$ over a given affine subspace. Usually, the subspace is described by a linear equations of the form $\mathscr{A}\mathscr{X}=\mathbf{b}$ as discussed in Section \ref{Sect. 5.1}. This problem can be formulated as a convex optimization in the following form:
\begin{eqnarray}\label{nuclear-norm-min}
\min\limits_{\mathscr{X}}\; \|\mathscr{X}\|_*\quad
\mbox{\rm s.t.}\quad \mathscr{A}\mathscr{X}=\mathbf{b}.
\end{eqnarray}
Then, by using the TSDP characterization of the nuclear norm given in (\ref{nuclear-norm-TSDP-dual}),
we can rewrite (\ref{nuclear-norm-min}) as
\begin{eqnarray*}
\min\limits_{\mathscr{X},\mathscr{W}_1,\mathscr{W}_2}\;
\frac{1}{2p}[Tr(\mathscr{W}_1)+Tr(\mathscr{W}_2)]
\quad
\mbox{\rm s.t.}\quad \left[\begin{array}{ll}
\mathscr{W}_1 & \mathscr{X}\\
\mathscr{X}^\top &\mathscr{W}_2
\end{array}\right]\succeq_\mathscr{T} \mathscr{O},\;\;\mathscr{A}\mathscr{X}=\mathbf{b}.
%
\end{eqnarray*}}

It was showed in \cite{LFCL-19} that the tensor robust principal component analysis problem can be transformed into the nuclear norm minimization problem with such tensor nuclear norm, which can be solved with the help of the theory of tensor decomposition. The above discussion demonstrates that the minimization of the above tensor nuclear norm of a third-order tensor can also be solved by dealing with the corresponding TSDP problem. As is known to us, a lot of practical problems are usually turned out to be low-rank models with third-order tensors, and most of them can be solved by the nuclear norm minimization problem with another related tensor nuclear norm defined in \cite{HKBH-13}, which has been shown to be widely applied in some practical problems, such as image processing, tensor principal component analysis, tensor completion, and so on. Then it is worthwhile to investigate these forms of convex relaxation by replacing the nuclear norm defined in \cite{HKBH-13} by the one given in \cite{LFCL-19}, and the TSDP provides another path to achieve solutions of low-rank recovery problems with third-order tensors appearing in the real-life applications.

{\bf Application 4. Integer quartic programming.} Analogue to the classic SDP relaxation of integer quadratic programming, we investigate the TSDP relaxation for the following integer quartic programming:
\begin{eqnarray}\label{IQP}
\max\limits_\mathscr{X}\; \langle\mathscr{X},\mathscr{A}\ast\mathscr{X}\rangle\quad
\mbox{\rm s.t.}\quad \mathbf{X}=\mathbf{x}\mathbf{x}^\top,\;\; x_i\in\{+1,-1\},\;\; \forall\ i\in[n],
\end{eqnarray}
where $\mathscr{A}\in S\mathbb{R}^{n\times n\times n}$ is given, $\mathbf{x}:=(x_i:i\in [n])^\top\in \mathbb{R}^n$, $\mathbf{X}\in \mathbb{R}^{n\times n}$ and $\mathscr{X}\in \mathbb{R}^{n\times 1\times n}$ is the corresponding tensor of the matrix $\mathbf{X}$.

Noting that for any $i\in[n]$, $x_i\in\{+1,-1\}$ if and only if $x_i^2=1$, if and only if $x_i^2(x_1^2+x_2^2+\cdots+x_n^2)=n$. Thus, problem (\ref{IQP}) is equivalent to
\begin{eqnarray}\label{IQP-2}
\quad \max\limits_\mathscr{X}\; \langle\mathscr{X},\mathscr{A}\ast\mathscr{X}\rangle\quad
\mbox{\rm s.t.}\quad \mathbf{X}=\mathbf{x}\mathbf{x}^\top,\;\; x_i^2(x_1^2+x_2^2+\cdots+x_n^2)=n,\;\; \forall\ i\in[n].
\end{eqnarray}

Now, we try to arrive at the TSDP relaxation of (\ref{IQP-2}) by using the Lagrangian multiplier method. By adding Lagrangian multiplier $y_i\in \mathbb{R}$ to each equality constraint in (\ref{IQP-2}), we can obtain the Lagrangian function for (\ref{IQP-2}):
$$\begin{array}{rcl}
L(\mathscr{X},\mathbf{y})&:=&\langle\mathscr{X},\mathscr{A}\ast\mathscr{X}\rangle-\sum\limits_{i=1}^n y_{i}(x_i^2(x_1^2+x_2^2+\cdots+x_n^2)-n)\\
&=&\langle\mathscr{X},\mathscr{A}\ast\mathscr{X}\rangle-\langle\mathscr{X},\mathscr{D}iag(\mathbf{y})\ast\mathscr{X}\rangle+n\mathbf{e}^\top\mathbf{y},
\end{array}$$
where $\mathscr{D}iag(\mathbf{y})$ represents the $F$-diagonal tensor induced by $\mathbf{y}:=(y_i:i\in [n])^\top$ with $\mathscr{D}iag(\mathbf{y})^{(1)}$ being the diagonal matrix $Diag(y_i:i\in [n])$ and other frontal slices being zeroes, and $\mathbf{e}:=(1,1,\ldots,1)^\top\in \mathbb{R}^n$. Since
$$\begin{array}{rcl}
\max\limits_\mathscr{X} L(\mathscr{X},\mathbf{y})&=&\max\limits_\mathscr{X}[\langle\mathscr{X},\mathscr{A}\ast\mathscr{X}\rangle-\langle\mathscr{X},\mathscr{D}iag(\mathbf{y})\ast\mathscr{X}\rangle+n\mathbf{e}^\top\mathbf{y}]\\
&=&-\min\limits_{\mathscr{X}}[\langle\mathscr{X},(\mathscr{D}iag(\mathbf{y})-\mathscr{A})\ast\mathscr{X}\rangle-n\mathbf{e}^\top\mathbf{y}]\\
&=& \left\{\begin{array}{ll}
n\mathbf{e}^\top\mathbf{y}, & \mbox{\rm if}\;\; \mathscr{D}iag(\mathbf{y})-\mathscr{A}\succeq_\mathscr{T}\mathscr{O},\\
 +\infty, & \mbox{\rm otherwise},
\end{array}\right.
\end{array}$$
the Lagrangian dual problem of (\ref{IQP-2}) turns out to be
\begin{eqnarray*}
\min\limits_{\mathbf{y}}\; n\mathbf{e}^\top\mathbf{y} \quad
\mbox{\rm s.t.}\quad \mathscr{D}iag(\mathbf{y})-\mathscr{A}\succeq_\mathscr{T}\mathscr{O},
\end{eqnarray*}
which is a TSDP model and its dual problem is
\begin{eqnarray*}
\max\limits_{\mathscr{X}}\; \mathscr{A}\bullet \mathscr{X} \quad
\mbox{\rm s.t.}\quad \mathscr{X}_{ii1}=n, \;\; i\in [n],\;\; \mathscr{X}\succeq_\mathscr{T} \mathscr{O}.
\end{eqnarray*}

Previously, we have converted some optimization problems over tensor space and matrix space into TSDP problems, respectively. Next, we show that some optimization problems over vector space can also be solved by TSDP problem models, such as polynomial optimization problems.

{\bf Application 5. Calculating the global lower bound of a polynomial of even degree.}
Consider the polynomial optimization problem:
\begin{eqnarray}\label{f}
f^{uc}:=\min_{\mathbf{x}\in \mathbb{R}^n}\; f(\mathbf{x})=f_0+\sum_{\alpha\in U^n_{2d}}{f_\alpha}\mathbf{x}^\alpha,
\end{eqnarray}
where $U^n_{2d}=\{\alpha\in \mathbb{N}^n: 0<|\alpha|\leq 2d\}$. As is well-known to us, problem (\ref{f}) can be solved by a relaxation into the following model through the sums of squares (SOS for short) method \cite{P-03}:
\begin{eqnarray}\label{f-SOS-SDP}
f^{uc}_{sos}:=\max\limits_{(\mathbf{x},\gamma)\in \mathbb{R}^n\times \mathbb{R}}\; {\gamma} \quad \mbox{\rm s.t.}\quad f(\mathbf{x})-\gamma \;\;\mbox{\rm is \ SOS}.
\end{eqnarray}
Actually, define $[\mathbf{x}]_d=(1,x_1,\ldots,x_n,x_1^2,x_1x_2,\ldots,x_1x_n,x_2^2,x_2x_3,\ldots,x_n^2,\ldots x_1^d,\ldots, x_n^d)^\top$,
which is a column vector of size $C_{n+d}^d$ consisting of all monomials whose degrees are no more than $d=\frac{deg(f)}{2}$, then (\ref{f-SOS-SDP}) can be transformed into a standard SDP as \cite{NW-12}:
\begin{eqnarray}\label{f-SDP}
f^{uc}_{sdp}:=\max\limits_{\mathbf{X}}\; \mathbf{C}\bullet\mathbf{X} \quad \mbox{\rm s.t.}\quad \mathbf{A}\mathbf{X}=\mathbf{b},\;\; {\mathbf{X} \succeq \mathbf{O}},
\end{eqnarray}
where $\mathbf{b}=({f_\alpha})_{\alpha\in U^n_{2d}}$ whose dimension is $C_{n+2d}^{2d}-1$, and $\mathbf{A}$ is a linear operator denoted as $\mathbf{A}\mathbf{X}=(\mathbf{A_\alpha}\bullet\mathbf{X})_{\alpha\in U^n_{2d}}$ with $\mathbf{A_\alpha}$ and $\mathbf{C}$ being constant symmetric matrices such that
$[\mathbf{x}]_d[\mathbf{x}]_d^\top=\mathbf{C}+\sum_{\alpha\in U^n_{2d}}\mathbf{A_\alpha}\mathbf{x}^\alpha$. In addition, $f^{uc}_{sos}=f_0-f^{uc}_{sdp}$.

As is known to us, in the process of solving the polynomial optimization problem with the SDP relaxation, one of the challenges is that the size of the SDP model increases significantly with the increasing of the number of variables or the degree of polynomial. Next, we show that the minimization problem of $f$ can be relaxed into a standard TSDP problem by rearranging the above monomial vector $[\mathbf{x}]_d$ into a third-order tensor form and then exploiting the properties of T-PSD tensors, which can be dealt with by solving an SDP problem of smaller size than the one in (\ref{f-SDP}).

Suppose that $C_{n+d}^d=mp$, then $[\mathbf{x}]_d\in \mathbb{R}^{mp}$. Let $[\mathscr{X}]_d$ be a tensor in $\mathbb{R}^{m\times 1\times p}$ with $[\mathscr{X}]_d=fold([\mathbf{x}]_d)$. Then by Theorem \ref{equivalent-PSD}, a tensor $\mathscr{A}\in S\mathbb{R}^{m\times m\times p}$ is T-PSD iff there exists some tensor $\mathscr{P}\in S\mathbb{R}^{l\times m\times p}$ such that $\mathscr{A}=\mathscr{P}^\top\ast \mathscr{P}$. Thus, $f(\mathbf{x})-\gamma$ must be SOS (and then be nonnegative) if there exists $\mathscr{X}\in S\mathbb{R}^{m\times m\times p}$ such that
\begin{eqnarray}\label{SOS-TSDP}
f(\mathbf{x})-\gamma=\frac{1}{p}Tr([\mathscr{X}]_d^\top\ast\mathscr{X}\ast [\mathscr{X}]_d)=\mathscr{X}\bullet([\mathscr{X}]_d\ast[\mathscr{X}]_d^\top),\quad \mathscr{X}\succeq_\mathscr{T}\mathscr{O}.
\end{eqnarray}
Therefore, finding the minimum value of $f$
can be relaxed into the following problem:
\begin{eqnarray}\label{f-SOS-TSDP}
f^{uc}_{sos-tsdp}:=\max\; {\gamma} \quad \mbox{\rm s.t.}\quad f(\mathbf{x})-\gamma=\mathscr{X}\bullet([\mathscr{X}]_d\ast[\mathscr{X}]_d^\top), \;\; \mathscr{X}\succeq_\mathscr{T}\mathscr{O}.
\end{eqnarray}
Define constant symmetric tensors $\mathscr{C}$ and $\mathscr{A_\alpha}$ such that
$$[\mathscr{X}]_d\ast[\mathscr{X}]_d^\top=\mathscr{C}+\sum_{\alpha\in U^n_{2d}}\mathscr{A_\alpha}\mathbf{x}^\alpha,$$
then (\ref{SOS-TSDP}) can be expressed as follows:
$$f(\mathbf{x})-\gamma=\mathscr{C}\bullet\mathscr{X}+\sum_{\alpha\in U^n_{2d}}(\mathscr{A_\alpha}\bullet\mathscr{X})\mathbf{x}^\alpha,  \quad \mathscr{X}\succeq_\mathscr{T}\mathscr{O}.$$
Noting that $f(\mathbf{x})=f_0+\sum_{\alpha\in U^n_{2d}}{f_\alpha}\mathbf{x}^\alpha$,
hence $\gamma$ is feasible for (\ref{f-SOS-TSDP}) if and only if there exists $\mathscr{X}\succeq_\mathscr{T}\mathscr{O}$ such that
$\mathscr{C}\bullet\mathscr{X}+\gamma = f_0$, and $\mathscr{A_\alpha}\bullet\mathscr{X} = f_\alpha$ for any $\alpha\in U^n_{2d}$.
Define a linear operator from $S\mathbb{R}^{m\times m\times p}$ into $\mathbb{R}^{(C_{n+2d}^{2d}-1)}$ as $\mathscr{A}\mathscr{X}=[\mathscr{A_\alpha}\bullet\mathscr{X}]_{\alpha\in U^n_{2d}}$. Then up to a constant, the problem (\ref{f-SOS-TSDP}) is equivalent to the TSDP problem:
\begin{eqnarray}\label{f-TSDP}
f^{uc}_{tsdp}:=\max\limits_{\mathscr{X}}\; \mathscr{C}\bullet\mathscr{X} \quad \mbox{\rm s.t.}\quad \mathscr{A}\mathscr{X}=\mathbf{b},\quad {\mathscr{X} \succeq_\mathscr{T} \mathscr{O}},
\end{eqnarray}
and $f^{uc}_{sos-tsdp}=f_0-f^{uc}_{tsdp}$.

Besides, from the above discuss, it is easy to see that (\ref{f-SOS-TSDP}) is both relaxation of (\ref{f}) and (\ref{f-SOS-SDP}) when $p\neq 1$, i.e., $f^{uc}\geq f^{uc}_{sos}\geq f^{uc}_{sos-tsdp}$, and (\ref{f-SOS-TSDP}) reduces to (\ref{f-SOS-SDP}) when $p=1$. Thereby, it should be noticed that (\ref{f-TSDP}) with $p\neq 1$ is a further relaxation of (\ref{f-SDP}). Then {\it whether or not is there a case where the relaxation of (\ref{f-TSDP}) with $p\neq 1$ can achieve the same effect as (\ref{f-SDP})?} Below, we answer this question by giving a theorem to show the necessary and sufficient condition of $f^{uc}_{sdp}=f^{uc}_{tsdp}$ under the assumption that $p\neq 1$.

\begin{theorem}\label{SDP=TSDP}
Suppose that $f(\mathbf{x}):\mathbb{R}^n\rightarrow R$ is a $2d$-degree polynomial, whose global lower bound can be solved by both an SDP relaxation as (\ref{f-SDP}) and a TSDP relaxation as (\ref{f-TSDP}) with $p\neq1$. Then $f^{uc}_{sdp}=f^{uc}_{tsdp}$ if and only if there exists an optimal solution $\mathbf{X^*}$ of (\ref{f-SDP}) which is $p$-block circular.
\end{theorem}

\begin{proof} Recall that (\ref{f-TSDP}) is equivalent to (\ref{f-SOS-TSDP}), that is
\begin{eqnarray}\label{1}
f_0-f^{uc}_{tsdp}:=\max\; {\gamma} \quad \mbox{\rm s.t.}\quad f(\mathbf{x})-\gamma=\mathscr{X}\bullet([\mathscr{X}]_d\ast[\mathscr{X}]_d^\top), \;\; \mathscr{X}\succeq_\mathscr{T}\mathscr{O};
\end{eqnarray}
and (\ref{f-SDP}) is equivalent to (\ref{f-SOS-SDP}), which is further equivalent to the following problem:
\begin{eqnarray}\label{2}
f_0-f^{uc}_{sdp}:=\max\; {\gamma} \quad \mbox{\rm s.t.}\quad f(\mathbf{x})-\gamma=\mathbf{X}\bullet([\mathbf{x}]_d \cdot[\mathbf{x}]_d^\top), \;\; \mathbf{X}\succeq\mathbf{O}.
\end{eqnarray}
From $$
\begin{array}{lcl}
\mathscr{X}\bullet([\mathscr{X}]_d\ast[\mathscr{X}]_d^\top)&=&\frac{1}{p}Tr([bcirc([\mathscr{X}]_d)]^\top\cdot bcirc(\mathscr{X})\cdot bcirc([\mathscr{X}]_d))\\
&=&\frac{1}{p}bcirc([\mathscr{X}]_d)\bullet [bcirc(\mathscr{X})\cdot bcirc([\mathscr{X}]_d)]\\
&=&(\frac{1}{p}\cdot p)unfold([\mathscr{X}]_d)\bullet [bcirc(\mathscr{X})\cdot unfold([\mathscr{X}]_d)]\\
&=&[\mathbf{x}]_d\bullet[ bcirc(\mathscr{X})\cdot [\mathbf{x}]_d]\\
&=& bcirc(\mathscr{X})\bullet([\mathbf{x}]_d\cdot[\mathbf{x}]_d^\top),
\end{array}$$
it is easy to see that if $\mathscr{X}^*$ is an optimal solution of (\ref{1}), then $bcirc(\mathscr{X}^*)$ must be a feasible solution of (\ref{2}) and $f^{uc}_{sdp}\leq f^{uc}_{tsdp}$.

On one hand, suppose that $\mathscr{X}^*$ is an optimal solution of (\ref{1}), if $f^{uc}_{sdp}=f^{uc}_{tsdp}$, then $bcirc(\mathscr{X}^*)$ is exact an optimal solution of (\ref{2}) and  (\ref{f-SDP}), which is a $p$-block circular matrix. On the other hand, if there exists an optimal solution $\mathbf{X^*}$ of (\ref{f-SDP}) which is $p$-block circular, then by choosing $\mathscr{X}^*=bcirc^{-1}(\mathbf{X^*})$ we can find that $\mathscr{X}^*$ is a feasible solution of (\ref{1}) and the corresponding value of $\gamma$ is $f_0-f^{uc}_{sdp}$ when $\mathscr{X}=\mathscr{X}^*$ in (\ref{1}). Noting that the optimal value of  (\ref{1}) is $f_0-f^{uc}_{tsdp}$, then we have that $f_0-f^{uc}_{sdp}\leq f_0-f^{uc}_{tsdp}$, which together with $f^{uc}_{sdp}\leq f^{uc}_{tsdp}$ implies that $f^{uc}_{sdp}= f^{uc}_{tsdp}$.
\end{proof}

Theorem \ref{SDP=TSDP} states that if there exists an optimal solution $\mathbf{X^*}$ of (\ref{f-SDP}) is $p$-block circular, then the relaxation problem (\ref{f-TSDP}) can achieve not worse effect than (\ref{f-SDP}). Below, we show that there are some benefits in computation costs and storage costs of (\ref{f-TSDP}) compared with (\ref{f-SDP}) for the cases where the unconstrained polynomial optimization can be solved by both TSDP relaxation as (\ref{f-TSDP}) and SDP relaxation as (\ref{f-SDP}):
%
\begin{itemize}
\item {\bf The number of entries of the variable $\mathscr{X}$ in (\ref{f-TSDP}) is less than the one of $\mathbf{X}$ in (\ref{f-SDP}).} Since $[\mathscr{X}]_d\in \mathbb{R}^{m\times 1\times p}$ and $[\mathbf{x}]_d\in \mathbb{R}^{mp}$, it is easy to find that $\mathbf{X}$ in (\ref{f-SDP}) is a matrix with $mp\times mp$ entries, while $\mathscr{X}$ in (\ref{f-SDP}) turns out to be a tensor in $S\mathbb{R}^{m\times m\times p}$ which contains only $m\times m\times p$ entries.
\item {\bf (\ref{f-TSDP}) can be solved by transformed into a CSDP problem with block diagonal sparse structure, which is not available for (\ref{f-SDP}) even if the variable of (\ref{f-SDP}) possess special structure since the constant matrice $\mathbf{A}_\alpha$ is fixed and do not possess the block diagonal structures.} Specifically speaking, if we adopt the method described in Section \ref{Sect. 5.2} to deal with (\ref{f-TSDP}), we only need to solve a CSDP with $\frac{p+1}{2}$ or $\frac{p+2}{2}$ blocks of size $m\times m$ under some equality constraints, while (\ref{f-SDP}) can be seen as dealing with an ordinary SDP with $1$ block of size $mp\times mp$ under the same number of equality constraints. Hence, for the minimization problem of a real polynomial of even degree, solving via the TSDP relaxation built the above would lead to lower computation costs than the corresponding SDP relaxation, which can also be seen from the numerical experiments in the next subsection.
\item {\bf The constant symmetric tensors $\mathscr{A}_\alpha$ for $\alpha\in U^n_{2d}$ in (\ref{f-TSDP}) spend no more storage cost than those matrices $\mathbf{A}_\alpha$ in (\ref{f-SDP}).} Since the number of entries in the tensor $[\mathscr{X}]_d\ast[\mathscr{X}]_d^\top$ is $\frac{1}{p}$ of those in the moment matrix $[\mathbf{x}]_d[\mathbf{x}]_d^\top$ and each entry of $[\mathscr{X}]_d\ast[\mathscr{X}]_d^\top$ is consisted of a linear combination of $p$ monomials whose degrees are no more than $2d$. Thus it seems that we need to take as the same storage as those for matrices $\mathbf{A}_\alpha$ to store these tensors $\mathscr{A}_\alpha$. However, noting that some entries of $[\mathscr{X}]_d\ast[\mathscr{X}]_d^\top$ may be consisted of a linear combination of $p$ same monomials whose degrees are no more than $2d$, therefore the constant symmetric tensors $\mathscr{A}_\alpha$ for $\alpha\in U^n_{2d}$ in (\ref{f-TSDP}) could sometimes save some storage than those matrices $\mathbf{A}_\alpha$ in (\ref{f-SDP}). In addition, the storage capacity does not change in the transformation of TSDPs into CSDPs as shown in Section \ref{Sect. 5.2}.
\end{itemize}

\subsection{Numerical computation}\label{Sect. 5.4} 
In this subsection, we report preliminary numerical results for solving TSDPs by the method shown in Section \ref{Sect. 5.2}. Taking Application 5 discussed in Section \ref{Sect. 5.3} for example, we consider two polynomial optimization problems and implement these problems in Matlab R2016a on our PC via SDPNAL$+$ \cite{SDPNAL+-19}. The computation is performed on a Dell Laptop with CPU of 3.2 GHz and RAM of 4.0 GB. As pointed in Application 5 of Section \ref{Sect. 5.2}, the TSDP relaxation of unconstrained polynomial optimization can achieve the same optimal value as the SDP relaxation sometimes, thus here we confirm this conclusion by two simple examples and illustrate the benefits of the TSDP relaxation furthermore.

\noindent {\bf Example 1.} Minimize the following polynomial:
$$\begin{array}{rcl}
f(\mathbf{x})&=&(x_1+x_2^3+x_1^2x_2)^2+(x_1+x_1^2+x_2^3)^2+(x_1+x_1^2+x_2^2)^2\\
&&+(x_1^2+x_2^2+x_1^2x_2)^2+(x_2^2+x_1^2x_2+x_2^3)^2.
\end{array}$$
It is obvious that the global minimum value $f^*=0$, the number of variables $n=2$, $d=\frac{deg(f)}{2}=3$ and $f(\mathbf{x})-f^*$ can be expressed as
$$f(\mathbf{x})-f^*=\left[\begin{array}{c}
                        1 \\
                        x_1 \\
                        x_2 \\
                        x_1^2 \\
                        x_1x_2 \\
                        x_2^2 \\
                        x_1^3 \\
                        x_1^2x_2 \\
                        x_1x_2^2 \\
                        x_2^3 \\
                      \end{array}\right]^\top
                      \left[\begin{array}{cc|cc|cc|cc|cc}
                                     0 & 0 & 0 & 0 & 0 & 0 & 0 & 0 & 0 & 0 \\
                                     0 & 3 & 0 & 2 & 0 & 1 & 0 & 1 & 0 & 2 \\ \hline
                                     0 & 0 & 0 & 0 & 0 & 0 & 0 & 0 & 0 & 0 \\
                                     0 & 2 & 0 & 3 & 0 & 2 & 0 & 1 & 0 & 1 \\ \hline
                                     0 & 0 & 0 & 0 & 0 & 0 & 0 & 0 & 0 & 0 \\
                                     0 & 1 & 0 & 2 & 0 & 3 & 0 & 2 & 0 & 1 \\ \hline
                                     0 & 0 & 0 & 0 & 0 & 0 & 0 & 0 & 0 & 0 \\
                                     0 & 1 & 0 & 1 & 0 & 2 & 0 & 3 & 0 & 2 \\ \hline
                                     0 & 0 & 0 & 0 & 0 & 0 & 0 & 0 & 0 & 0 \\
                                     0 & 2 & 0 & 1 & 0 & 1 & 0 & 2 & 0 & 3 \\
                                   \end{array}\right]
                      \left[\begin{array}{c}
                          1 \\
                        x_1 \\
                        x_2 \\
                        x_1^2 \\
                        x_1x_2 \\
                        x_2^2 \\
                        x_1^3 \\
                        x_1^2x_2 \\
                        x_1x_2^2 \\
                        x_2^3 \\
                      \end{array}\right],
$$
where the square matrix is $5$-block circular. Therefore, by Theorem \ref{SDP=TSDP}, we know that the TSDP relaxation with $p=5$ can achieve not worse effect than the SDP relaxation. Specifically, to solve this problem via the TSDP relaxation as described from (\ref{f}) to (\ref{f-TSDP}), we first find the corresponding $\mathscr{A}_\alpha$, $\mathscr{C}$ and $\mathbf{b}$ in (\ref{f-TSDP}), then transform $\mathscr{A}_\alpha$ and $\mathscr{C}$ into the corresponding $\mathbf{A}_\alpha$ and $\mathbf{C}$ as the procedure from (PTSDP) to (P$''$CSDP), and finally, call the package SDPNAL$+$ to solve the derived SDP. In addition, we also solve this problem via the traditional SDP relaxation as described from (\ref{f}) to (\ref{f-SDP}) for comparison.

Numerical results for {\bf Example 1} by SDP and TSDP via SDPNAL$+$ are shown in \ref{Table 1}, where the size of the TSDP is replaced by the size of the corresponding CSDP with block structure and {\bf in each pair ``($blk,N,m$)", $blk$, $N$ and $m$ mean} the number of blocks, the length of each block and the number of equality constraints, respectively; {\bf``opt" means} the computed optimal value of the corresponding relaxation; {\bf``cpu1" means} the time in seconds spent for finding the corresponding matrices $\mathbf{A}_\alpha$ of the SDP or the CSDP; {\bf``cpu2" means} the time in seconds spent for solving corresponding SDP problem or CSDP problem via SDPNAL$+$; and {\bf``cpu3" means} the total time in seconds spent for solving Example 1 by the SDP relaxation or the TSDP relaxation.

\begin{table}[!htbp]
  \centering{\small
     \begin{tabular}{ccccccccc}
     \hline
          & n & d  & p  & ($blk,N,m$)  &   opt                & cpu1(s)   &cpu2(s)   & cpu3(s)   \\
          \hline
     SDP  & 2 & 3 & 1  &  (1,10,28)   & $2.1545e^{-09}$       &0.021      & 0.969     &1.008   \\
     TSDP & 2 & 3 & 5 &   (5,2,28)    & $-1.1183e^{-15}$      &0.037      & 0.426     &0.483   \\
     \hline
\end{tabular}}
\caption{Numerical results for Example 1 by the SDP and the TSDP via SDPNAL$+$.}
\label{Table 1}
\end{table}

From \ref{Table 1}, we can see that the TSDP relaxation performs well for \textbf{Example 1}. Especially, for the TSDP relaxation, we obtain that the global minimum value $f^*_{tsdp}=-1.1183e^{-15}$ by solving a CSDP problem with block diagonal structure, which only takes less than one-twice of the time of the corresponding SDP relaxation. In addition, the TSDP relaxation has higher accuracy than the SDP relaxation in this example.

\noindent {\bf Example 2.} Minimize the following polynomial:
$$\begin{array}{ll}
f(\mathbf{x})=& 1+x_1^{10}x_2^4+x_1^{8}x_2^{12}+x_1^{24}x_2^{2}+x_1^{24}x_2^{6}+x_1^{32}x_2^{2}+x_1^{8}x_2^{28}+
x_1^{28}x_2^{12}+x_1^{10}x_2^{32}\\
&+x_1^{42}x_2^{4}+x_1^{30}x_2^{18}+x_1^{20}x_2^{30}+x_1^{12}x_2^{40}+x_1^6x_2^{48}+x_1^2x_2^{54}+x_2^{58}.
\end{array}$$
It is obvious that the global minimum value $f^*=1$, the number of variables $n=2$, and $d=\frac{deg(f)}{2}=29$. We test $p=3$ and $p=15$ for the TSDP relaxation of this problem, respectively. Especially, by the TSDP relaxation with $p=15$, we obtain that the global minimum value $f^*_{tsdp}=1+6.1507e^{-08}$ by solving a CSDP problem with block diagonal structure, which takes only 2.5 seconds; while the corresponding SDP relaxation takes 90.3 seconds. Meanwhile, the TSDP relaxation has higher accuracy than the SDP relaxation in this example.

Numerical results for \textbf{Example 2} by the SDP relaxation and the TSDP relaxation are shown in \ref{Table 2}, where ``($blk,N,m$)'', ``opt'' , ``cpu1'',``cpu2'' and ``cpu3'' are same as those in \ref{Table 1}.
\begin{table}[!htbp]
  \centering{\small
     \begin{tabular}{ccccccccc}
     \hline
          & n & d  & p   &  ($blk,N,m$)          &   opt         &cpu1(s)   & cpu2(s)  &cpu3(s) \\
          \hline
     SDP  & 2 & 29 & 1   & (1,465,1769)    &  $1-1.1897e^{-07}$  &109.453   & 90.303   &199.896 \\
     TSDP & 2 & 29 & 3   & (2,155,1769)    &  $1-4.3220e^{-08}$  &107.834   & 61.426   &170.272 \\
     TSDP & 2 & 29 & 15  & (8,31,1769)     &  $1+6.1507e^{-08}$  &103.302   & 2.515    &110.483 \\
     \hline
\end{tabular}}
\caption{Numerical results for Example 2 by the SDP and the TSDP via SDPNAL$+$.}\label{Table 2}
\end{table}

From Tables \ref{Table 1} and \ref{Table 2}, we can see that the TSDP relaxation has better performances both in time cost and precision than the traditional SDP relaxation for these test examples. Especially in the time cost, solving the CSDP problem transformed from the corresponding TSDP problem saves a lot of time compared with solving the corresponding SDP problem; while the time for finding the corresponding complex matrices $\mathbf{A}_\alpha$ in the CSDP problem is almost the same as the time for finding the corresponding matrices $\mathbf{A}_\alpha$ in the SDP problem. In fact, it should be noticed that the time for finding the corresponding matrices $\mathbf{A}_\alpha$ can be removed from the total time of solving above polynomial optimization problems via the TSDP or the SDP relaxation because these $\mathbf{A}_\alpha$ are fixed as long as $n$, $d$ and $p$ are given, and we show them in Tables \ref{Table 1} and \ref{Table 2} just to make the time spent for the TSDP relaxation and the SDP relaxation more clear.

\section{Concluding remarks}\label{Sect. 5}
In this paper, we aimed to generalize the SDP problem to the third-order tensor case. For this purpose, we first introduced the T-positive semidefiniteness of third-order symmetric tensors from the second-order discrimination condition of the convexity of the multi-vector real-valued function, and then, extended some of useful characterizations and properties of the symmetric PSD matrix to the third-order symmetric T-PSD tensor. After that, we replaced the variable in the classic SDP by the third-order symmetric tensor to introduce the TSDP, which was dealt with by converting it to a CSDP with the block structure. Finally, several examples with respect to the TSDP problem were shown and some numerical results of minimizing two polynomials via the TSDP relaxation were reported, which demonstrated that our method by the TSDP relaxation performs better than the traditional SDP relaxation for the test examples.

Some issues need to be studied in the future.
\begin{itemize}
\item[(i)] In Section \ref{Sect. 5.4}, we have just done some preliminary numerical experiments to test the feasibility and effectiveness of solving TSDPs by dealing with the corresponding CSDPs. Surprisingly, we find that for some unconstrained polynomial optimization, the method of the TSDP relaxation has good performance {sometimes}. Then what kind of polynomial optimization problems does the TSDP relaxation work well for? It deserves further study.
\item[(ii)] It is known to us that the SDP has shown great power in a very wide range of areas. In Section \ref{Sect. 5.3}, we have just presented a few simple transformations from some other models into TSDPs. We believe that more problems can be modeled (or relaxed) as TSDPs. It is also known that the T-product between third-order tensors promotes the emergence of many algorithms with good performance in many practical problems. It is possible that more efficient algorithms for TSDP problems can be designed by making use of the characteristics of the T-product and special-structures of practical problems. Furthermore, it deserves to study how to design efficient algorithms for solving large-scale  realistic problems.
\end{itemize}

\end{document}